\documentclass[a4paper,12pt]{article}

\begin{document}
\author{Sergey V. Ludkovsky.}
\title{Topological homeomorphisms of Banach spaces over non-Archimedean fields.}
\date{21 September 2013}
\maketitle
\begin{abstract}
The article is devoted to topological homeomorphisms of Banach
spaces over complete non-Archimedean normed infinite fields.
\footnote{ 
\par keywords: Banach space, field, homeomorphism  }
\end{abstract}

\par address: Department of Applied Mathematics,
\par Moscow State Technical University MIREA,
av. Vernadsky 78,
\par Moscow, 119454, Russia
\par e-mail: sludkowski@mail.ru

\section{Introduction.}
The topological classification of real Banach spaces was first made
by Bessaga and Pelczynski \cite{bess65,besspelcz62,besspelczb}. But
their methods and constructions were rather abstract. Later on
Anderson, Bing and Kadets have simplified their methods and have
constructed more concrete topological homeomorphisms of real Banach
spaces \cite{anderbing,kadets}.
\par Analogous situation is with Banach spaces over other fields with non-Archimedean
norms. It can be demonstrated that they are homeomorphic with
products of discrete topological spaces of suitable cardinality.
Although they can be classified by methods of the general topology,
but they are too abstract to be useful for applications in
non-Archimedean functional analysis, measure theory, topological
algebra or geometry. Indeed, all non topological structures of
functional analysis, measure theory, topological algebra or geometry
are lost in general topology. On the other hand, structures of
fields, even local, are very important in the aforementioned
branches of mathematics.
\par This article is devoted to investigations and constructions of locally analytic
homeomorphisms of Banach spaces over fields with non-Archimedean
norms. Norms are considered not only with values in the real field,
but also in linearly ordered rings.
\par Mathematical analysis over infinite non-Archimedean normed fields is
developing fast during recent years \cite{roo,sch1,vla3,yojan}.
Certainly structure of totally disconnected topological spaces, as
well as manifolds and Banach spaces over non-Archimedean infinite
fields is important not only for general topology, but also for
mathematical analysis. Previously non-archimedean polyhedral
expansions of totally disconnected $T_1\cap T_{3.5}$ topological
spaces and non-archimedean Banach manifolds were studied in
\cite{luumpe,lufpmpe,luum985,ludqimb}.
\par The Banach space $c_0$ in the non-archimedean analysis plays
the same principal role as the Hilbert space $l_2$ in the classical
analysis over $\bf R$ (see Theorems 5.13 and 5.16 in \cite{roo}).
\par In this paper locally analytic homeomorphisms of Banach
spaces over complete non-Archimedean normed infinite fields with
products of copies of the fields are investigated (see Theorem 2
below). \par All main results of this paper are obtained for the
first time.

\section{Homeomorphisms of linear topological spaces.}
\par {\bf 1. Notation.} For a topological space $X$ supplied with
a metric $\rho $ defining its topology let $B(X,y,r)$ denote a ball
$\{ z\in X: ~ \rho (z,y)\le r \} $ of radius $r>0$ containing a
marked point $y$. Let $\bf K$ be a field with a non-Archimedean
norm. This means that it is supplied with a multiplicative norm
satisfying the non-Archimedean inequality
\par $(1)$ $|x+y|\le \max (|x|,|y|)$ for all elements $x, y\in {\bf K}$.
\par Norms with values in $[0,\infty )$ are considered below, where as
usually $[0,\infty ) := \{ x: x\in {\bf R}, 0\le x \} .$ Such norm
is called non-Archimedean. With a non-Archimedean norm is associated
a multiplicative group $\Gamma _{\bf K} := \{ |x|: ~ x \in {\bf K},
x\ne 0 \} $. \par For example, ${\bf Q}_p\subset {\bf K}$ or ${\bf
F}_p(\theta ) \subset {\bf K}$, where ${\bf Q}_p$ denotes the field
of $p$-adic numbers, ${\bf F}_p(\theta )$ denotes the field of
formal Laurent series over the finite field ${\bf F}_p$, where $p>1$
is a prime number. Each $p$-adic number $x\in {\bf Q}_p$ has the
decomposition $x=\sum_{n=N}^{\infty } a_n p^n$, where $a_n\in \{ 0,
1,...,p-1 \} $ for each $n$, $N=N(x)\in {\bf Z}$ is an integer
number so that $a_N\ne 0$, $|x|=|x|_p =p^{-cN}$, $c>0$ is a constant
on ${\bf Q}_p$, particularly, $c=1$ can be taken, as usually ${\bf
Z}$ denotes the ring of all integers. The characteristic of the
$p$-adic field is zero $char ({\bf Q}_p)=0$, while the
characteristic of ${\bf F}_p(\theta )$ is $p= char ({\bf F}_p(\theta
) ) >1$ positive.
\par Each element $x\in {\bf F}_p(\theta )$ has the series
decomposition $x=\sum_{n=N}^{\infty } b_n \theta ^n$, where
$N=N(x)\in {\bf Z}$, $b_n\in {\bf F}_p$, $b_N\ne 0$, $|x| = |x|_p =
p^{-cN}$, $c>0$ is a constant on the field ${\bf F}_p(\theta )$,
particularly, $c=1$ can be taken. These fields differ by their
multiplication and addition rules. \par Then it is possible to take
an algebraic or a transcendental extension of an initial field and
its uniform or spherical completion if it is not such. The field
${\bf C}_p$ of complex $p$-adic numbers is obtained as the uniform
completion of a field containing all finite algebraic extensions of
the $p$-adic field ${\bf Q}_p$ with the norm extending that of on
${\bf Q}_p$. The fields ${\bf Q}_p$ and ${\bf F}_p(\theta )$ are
locally compact, the field ${\bf C}_p$ is not locally compact. The
normalization group $\Gamma _{\bf K} := \{ |x|: ~ x\in {\bf K}, x\ne
0 \} $ is multiplicative and commutative, for ${\bf K}={\bf Q}_p$ it
is isomorphic with $\{ p^n: ~ n\in {\bf Z} \} $, for ${\bf K}={\bf
C}_p$ it is isomorphic with $\{ p^x: ~ x\in {\bf Q} \} $, where
${\bf Q}$ denotes the field of all rational numbers. A larger field
${\bf U}_p$ exists so that ${\bf C}_p$ can be isometrically embedded
into ${\bf U}_p$ and the normalization group $\Gamma _{{\bf U}_p}$
is isomorphic with $ \{ p^x: ~ x\in {\bf R} \} $ (see
\cite{diar,sch1,weil}).
\par More generally a field $\bf K$ having a multiplicative
norm $|x|$ with values in a linearly ordered commutative topological
ring $\cal R$ can be considered so that $|x+y|\le \max (|x|, |y|)$.
Suppose that the ring $\cal R$ is complete as the uniform space (see
\cite{eng}) and $0$ denotes the neutral element relative to the
addition and $1$ is the unit element relative to the multiplication
in $\cal R$; moreover, if $p>1$ in ${\cal R}$, then $\lim_{n\to
+\infty }p^{-n}=0$, where $n\in {\bf N}$. We consider the case, when
an element $x\in {\bf K}$ exists so that $|x|=p>1$. For example,
${\cal R}\subset {\bf R}^{\gamma }$, where $\gamma $ is an ordinal,
while elements $z$ in ${\bf R}^{\gamma }$ are ordered
lexicographically: $y< z$, if $y_j=z_j$ for each $j<k$ and $y_k<z_k$
for some $k\in \gamma $, where ${\bf R}$ denotes the field of real
numbers. Non-zero elements $x\ne 0$ in the field ${\bf K}$ have
norms belonging to the multiplicative group $G$ of the ring $\cal
R$. We suppose that the normalization group $\Gamma _{\bf K} := \{
|x|: ~ x\in {\bf K}\setminus \{ 0 \} \} $ of the field $\bf K$ is
infinite and the closure of $\Gamma _{\bf K}$ in the ring $\cal R$
contains the zero point $0$, naturally, $\Gamma _{\bf K}\subset \{
z\in {\cal R}: ~ z>0 \} $; also $|x|=0$ if and only if $x=0\in {\bf
K}$.
\par Traditionally $c_0(\alpha
,{\bf K}) =: c_0(\alpha )$ denotes the normed space over the field
${\bf K}$ consisting of all nets $x=\{ x_j: ~ j\in {\alpha }, ~
x_j\in {\bf K} \} $ so that for each $\epsilon >0$ the set $ \lambda
(x,\epsilon ) := \{ j\in \alpha : ~ |x_j|>\epsilon \} $ is finite,
where $\alpha $ is a set, the norm of $x$ is:
\par $(2)$ $\| x \| := \sup_{j\in \alpha } |x_j|$. \\
That is, either $|\alpha |<\aleph _0$ or $\Gamma _{\bf K}$ is
discrete in ${\bf R}$ or $\cal R$ relative to the interval topology
induced by its linear ordering.
\par For a normed $\bf K$-linear space $E$
two vectors $x, y\in E$ are called orthogonal, if $\| ax+by \| =
\max (\| ax \| , \| by \| )$ for all $a, b\in \bf K$. For a real
number $0<t\le 1$ a finite or an infinite sequence of elements
$x_j\in E$ is called $t$-orthogonal, if $\| a_1x_1+....+a_mx_m +...
\| \ge t \max (\| a_1x_1 \| ,..., \| a_mx_m \|,... )$ for each
$a_1,...,a_m,... \in \bf K$ with $a_1x_1+...+a_mx_m+...\in E$.
\par The standard orthonormal in the non-Archimedean sense base in
$c_0(\alpha ,{\bf K})$ is $ e_j := (0,...,0,1,0,...)$ with $1$ at
the $j$-th place. The space $c_0(\alpha ,{\bf K})$ is Banach, when a
non-Archimedean field $\bf K$ is complete as a uniform space.
Henceforward, we consider the field $\bf K$ complete as the uniform
space, if something other will not be specified.
\par Let $\omega _0$ denote the first countable ordinal, for example,
${\bf N} := \{ 1,2,3,... \} $.
\par We consider the product $\prod_{j\in \alpha } X_j$ of
topological spaces $X_j$ supplied with the Tychonoff (product)
topology $\tau _{ty}$ with the base $U=\prod_{j\in \alpha }U_j$,
where each $U_j$ is open in $X_j$ and only a finite number of $U_j$
is different from $X_j$ for a given $U$ (see \cite{eng}).
Henceforth, a locally compact field $\bf K$ is considered so that
$\Gamma _{\bf K}$ is discrete. In this case let $p$ be such that
\par $p^{-1} := \sup \{ |x|: ~ x \in {\bf K}, ~ |x|< 1 \} $.

\par {\bf 2. Theorem.} {\it The Banach space
$c_0(\omega _0,{\bf K}) =: c_0$ over a uniformly complete
non-Archimedean field ${\bf K}$ supplied with its norm topology
$\tau _n$ is topologically homeomorphic with the countable product
${\bf K}^{\omega _0}$ of the non-Archimedean normed infinite field
${\bf K}$, where ${\bf K}^{\omega _0}$ is supplied with the
Tychonoff topology $\tau _{ty}$.}

\par The proof of this theorem is divided into several steps.
\par {\bf 3. Remark.} The condition that the field is infinite
is essential. If the field is finite, then it is discrete and
compact, consequently, the product ${\bf K}^{\omega _0}$ of compact
topological spaces is compact. But the linear topological space
${\bf c}_0(\omega _0,{\bf K})$ is not compact, since the covering by
clopen (closed and open simultaneously) balls $B({\bf c}_0(\omega
_0,{\bf K}),e_jx_m,1/p)$, $j\in \omega _0$, is infinite and has not
any finite sub-covering, where $\{ x_m: m \} $ are all distinct
elements of the field $\bf K$, $|x|=1$ for each $x\ne 0$, $|0|=0$,
$p>1$.
\par On the other hand, for the cardinality
$card (\alpha )> \aleph _0 := card (\omega _0)$ of the set $\alpha $
greater than $\aleph _0$ the base of neighborhoods of zero in ${\bf
K}^{\alpha }$ is uncountable, but $c_0(\alpha ,{\bf K})$ has a
countable base of neighborhoods of zero. Therefore, the topological
spaces $c_0(\alpha ,{\bf K})$ and ${\bf K}^{\alpha }$ are not
homeomorphic when $card (\alpha )> \aleph _0$.
\par A topology on the product $\prod_{j\in \alpha } X_j$
of spaces stronger than the Tychonoff topology is given by the base
$U = \prod_{j\in \alpha } U_j$, where each $U_j$ is open in $X_j$.
This topology is called the box topology $\tau _b$ \cite{nari}.
\par We use the notation ${\bf s} := {\bf K}^{\omega _0}$, ${\bf s}^{\alpha } :=
\prod_{j\in \alpha } {\bf K}_j$ for a subset $\alpha \subset \omega
_0$, where $\bf s$ and ${\bf s}^{\alpha }$ are supplied with the
product Tychonoff topology.

\par {\bf 4. Lemma.} {\it If $\alpha \subset \omega _0 $
and $ \beta = \omega _0 \setminus \alpha $ are disjoint subsets in
$\omega _0$ and $\alpha \ne \emptyset $ is non-void, then
\par $(1)$ $\bf s$ and ${\bf s}^{\alpha }\times {\bf s}^{\beta }$ are homeomorphic;
\par $(2)$ $c_0(\omega _0,{\bf K})$ and $c_0({\alpha },{\bf K})\times c_0({\beta },{\bf K})$
are homeomorphic. \par $(3)$. Moreover, if $\alpha $ is infinite,
then ${\bf s}^{\alpha }$ is homeomorphic with $\bf s$, while
$c_0({\alpha },{\bf K})$ is homeomorphic with $c_0(\omega _0,{\bf
K})$.}
\par $(1,3)$. This Lemma is evident, since $card (\omega _0^2)=card (\omega _0)=
\aleph _0$. \par $(2)$. If $x\in c_0(\alpha ,{\bf K})$ and $y\in
c_0(\beta ,{\bf K})$, then $\lim_j x_j=0$ when $\alpha $ is infinite
and $\lim_k y_k=0$ when $\beta $ is infinite, hence $\lim_l z_l=0$,
where $z_l=x_l$ for $l\in \alpha $, $ ~ z_l=y_l$ for $l\in \beta $,
also $\alpha \cup \beta = \omega _0$. At the same time $\| z \| =
\sup_{l\in \omega _0} |z_l| = \max ( \| x \| , \| y \| )$.

\par {\bf 5. Definitions.} Let $\pi _{\alpha }: c_0(\omega _0,{\bf K})\to
c_0(\alpha ,{\bf K})$ or $\pi _{\alpha }: {\bf s}\to {\bf s}^{\alpha
}$ denote the natural projection. In particular, if $\alpha = \{ j
\} $ is a singleton we can write $\pi _j$ instead of $\pi _{\alpha
}$.
\par A subset $E$ in $c_0(\alpha ,{\bf K})$ or in ${\bf s}^{\alpha }$ is called
deficient in the $j$-th direction, if $\pi _j(E)$ is a singleton
(i.e. consists of a single element). A subset $E$ in $c_0(\alpha
,{\bf K})$ or in ${\bf s}^{\alpha }$ is called infinitely deficient
if for some infinite subset $\beta \subset \alpha $, each projection
$\pi _j(E)$ is a singleton for every $j\in \beta $. In such case $E$
will also be called deficient with respect to $\beta $. Henceforth,
the term mapping will be used for continuous functions.
\par Suppose that $A$ is a topological space and $B$ is
its subset and $f$ is a mapping of $B$ into $A$. One says that the
mapping $f$ is limited by an open covering $W$ of $A$ if for each
point $x\in B$ there exists an element $V_x\in W$ so that $x\in V_x$
and $f(x)\in V_x$. \par If $A$ is a metric space supplied with a
metric $\rho $, then the supremum $\sup_{V\in W} diam (V)$ is called
the mesh of a covering $W$, where $diam (V) :=\sup_{a, b \in V} \rho
(a,b)$.
\par Let $f_1, f_2,...,f_m,...$ be a sequence of mappings such that
the limit $\lim_{m\to \infty } f_m\circ f_{m-1}\circ ... f_1: X\to
X$ exists, where $X$ is a topological space. This limit is denoted
by $L\prod_{j=1}^{\infty } f_j$ and is called the infinite left
product of the mappings $f_j$.
\par We consider the subsets $A_0 := \{ x\in c_0: ~
\sup_{i\in {\bf N}} |\sum_{j=1}^i x_j| =1 \} $ and $E^j := \{ x\in
c_0: ~ x=(x_1,x_2,...), ~ x_k=0 ~ \forall k>j \} $ in $c_0:=
c_0(\omega _0,{\bf K})$.
\par {\bf Lemma 6.} {\it The topological spaces ${\bf K}$ and
$B({\bf K},0,1)\setminus \{ 1 \} $ are topologically homeomorphic,
where ${\bf K}$ is a field (see \S 1).}
\par {\bf Proof.} We take any element $p\in G$ so that $p>1$ and
$p=|x|$ for some invertible element $x\in {\bf K}$ (see \S 1).
Therefore, $1/p^n=p^{-n}$ tends to zero in the topological ring
${\cal R}$ while $n\in {\bf N}$ tends to $+ \infty $. Then the
topological field ${\bf K}$ can be written as the disjoint union of
clopen subsets $B({\bf K},0,1)$, $B({\bf K},0,p)\setminus B({\bf
K},0,1)$,...,$B({\bf K},0,p^{n+1})\setminus B({\bf K},0,p^n),...$
with $n\in {\bf N}$. \par The norm in ${\bf K}$ is multiplicative,
consequently, $B({\bf K},0,p^n) = x^n B({\bf K},0,1)$ for each $n\in
{\bf Z}$, where $XY := \{ z=xy: ~ x\in X, y\in Y \} $ for two
subsets $X$ and $Y$ in ${\bf K}$. Thus subsets $B({\bf
K},0,p^{n+1})$ and $B({\bf K},0,p^m)$ are isomorphic for all $n,
m\in {\bf Z}$. Moreover, we have the equalities $B({\bf
K},0,1)\setminus \{ 1 \} = \bigcup_{n=0}^{\infty } [B({\bf
K},1,p^{-n})\setminus B({\bf K},1,p^{-n-1})]$ and $B({\bf
K},y,r)=y+B({\bf K},0,r)$ for each $y\in {\bf K}$ and $r>0$. Each
set $B({\bf K},1,p^{-n})\setminus B({\bf K},1,p^{-n-1})$ or $B({\bf
K},0,p^{n+1})\setminus B({\bf K},0,p^n)$ is the disjoint union of
balls $B({\bf K},y_j,p^m)$ or $B({\bf K},z_j,p^m)$ with $m= -n-1$ or
$m= n$ respectively. \par On the other hand, the quotient ring
$B({\bf K},0,1)/B({\bf K},0,1/p)$ exists. Its additive group is
isomorphic with \par $[x^nB({\bf K},0,1)]/[x^nB({\bf K},0,1/p)] =
B({\bf K},0,p^n)/B({\bf K},0,p^{n-1})$ \par $=x^n[B({\bf
K},0,1)/B({\bf K},0,1/p)]$ \\ for each $n\in {\bf Z}$, where $x\in
\bf K$ with $|x|=p$. Thus ${\bf K}$ and $B({\bf K},0,1)\setminus \{
1 \} $ can be presented as disjoint unions ${\bf K}=\bigcup_{\lambda
\in \Lambda _1}B({\bf K},z_{\lambda },r_{\lambda })$ and $B({\bf
K},0,1)\setminus \{ 1 \} =\bigcup_{\mu \in \Lambda _2} B({\bf
K},y_{\mu },r_{\mu })$ with $card (\Lambda _1)= card (\Lambda _2)
\ge \aleph _0$. We take any mapping $\phi : {\bf K}\to B({\bf
K},0,1)\setminus \{ 1 \}$ such that $\phi : B({\bf K},z_{\lambda
},r_{\lambda })\to B({\bf K},y_{\psi (\lambda )},r_{\psi (\lambda
)})$ is a homeomorphism. For example, $\phi $ can be chosen affine
$x\mapsto a+bx$ on each ball $B({\bf K},z_{\lambda },r_{\lambda })$,
where $\psi : \Lambda _1\to \Lambda _2$ is a bijective surjective
mapping. Thus $\phi : {\bf K}\to B({\bf K},0,1)\setminus \{ 1 \}$ is
the topological homeomorphism.
\par {\bf 7. Remark.} Using the preceding lemma we henceforth consider
${\bf s}$ as homeomorphic with $$(1)\quad {\bf s} \simeq
s=\prod_{j=1}^{\infty } [B({\bf K},0,1)_j\setminus \{ 1 \} ] $$
supplied with the Tychonoff product topology if $\bf s$ from \S 3
will not specified, where $B({\bf K},0,1)_j= B({\bf K},0,1)$ for
each $j$. Then we put $$(2)\quad s_* = \{ y\in s: ~ \lim_{m\to
\infty } (1-y_m)...(1-y_1)=0 \} \setminus \bigcup_{n=1}^{\infty }
{\sf E}^n,$$  where ${\sf E}^n := \{ x\in s: ~ x=(x_1,x_2,...), ~
x_k=0 ~ \forall k>n \} $.
\par {\bf 8. Lemma.} {\it A ring or a field ${\bf K}$ has the natural
uniformity and its completion $\tilde {\bf K}$ relative to this
uniformity is a topological ring or a field correspondingly.}
\par {\bf Proof.} The norm $|*|$ in $\bf K$ is multiplicative
with values in $G\cup \{ 0 \} \subset {\cal R}$. One can take a
diagonal $\Delta := \{ (x,y)\in {\bf K}^2: ~ x=y \} $ in the
Cartesian product ${\bf K}^2={\bf K}\times {\bf K}$. This norm
induces entourages of the diagonal in ${\bf K}^2$: $V_z := \{ (x,y)
\in {\bf K}^2: ~ |x-y|\le z \} $ for each $z\in G$. Therefore,
\par $(E1)$ $\bigcap_{z\in G} V_z=\Delta $, \\ since $x=y$ if and only if
$|x-y|=0$. \par $(E2).$ If $z_1<z_2$ then $V_{z_1}\subset V_{z_2}$,
since $|x-y|<z_1$ implies $|x-y|<z_2$. Then we have also
\par $(E3)$ if $(x,y)\in V_z$ and $(y,\xi )\in V_b$, then $(x,\xi )\in
V_{\max (z,b)}$, since $|x-\xi |\le \max (|x-y|,|y-\xi |)$.
\par Naturally, the inclusion $(x,y)\in V_z$ is equivalent to
$(y,x)\in V_z$, since $|x-y|=|y-x|$. The family $\cal E$ of all
entourages of the diagonal $D$ in $\bf K$ provides the uniformity
$\cal U$ in ${\bf K}$ compatible with its topology (see Chapter 8
\cite{eng}). The completion $\tilde {\bf K}$ relative to this
uniformity $\cal U$ is the uniformly complete field, since the
addition and multiplication operations are uniformly continuous on
the ring, also the inversion operation on ${\bf K}\setminus \{ 0 \}$
for the field and they have uniformly continuous extensions on
either $\tilde {\bf K}$ or on ${\tilde {\bf K}}\setminus \{ 0 \} $
respectively.
\par {\bf 9. Notation.} Let $\bf K$ be a uniformly complete
non-Archimedean field.  We define the subset \par $(1)$ $A_1 := \{
x\in c_0: ~ \sup_{k\in {\bf N}} |\sum_{j=1}^kx_j| =1,
|1-\sum_{j=1}^{k+1}x_j| \le |1-\sum_{j=1}^k x_j| ~ \forall k, ~
\sum_{j=1}^{\infty }x_j=1 \} $ in $c_0$, also \par $(2)$ $A_1^* :=
A_1\setminus \bigcup_{n=1}^{\infty } E^n$ (see \S 5).
\par Another larger subsets we define by the formula: \par  $(3)$ $A_2 :=
\{ x\in c_0: ~ \sup_{k\in {\bf N}} |\sum_{j=1}^kx_j| =1,
\sum_{j=1}^kx_j\ne 1 ~ \forall k, ~ \sum_{j=1}^{\infty }x_j=1 \} $
in $c_0$, also
\par $(4)$ $A_2^* := A_2\setminus \bigcup_{n=1}^{\infty } E^n$. \par Let
$A_1$ and $A_1^*$ and also $A_2$ and $A_2^*$ be supplied with the
topology inherited from the normed space $c_0$.
\par {\bf 10. Lemma.} {\it
The topological spaces $A_1^*$ and $s_*$ are homeomorphic.}
\par {\bf Proof.}  We define the following mapping $q: A_1^*\to s$ so
that $q(x)=y$, where $x=(x_1,x_2,...)\in A_1^*$, $y=(y_1,y_2,...)\in
s$ (see \S 7). The domain of $y_j$ is $B({\bf K},0,1)\setminus \{ 1
\} $ for each $j\in {\bf N}$.
\par If $x\in A_1$ and $|1-\sum_{j=1}^k
x_j|=0$, then $|1-\sum_{j=1}^mx_j|\le |1-\sum_{j=1}^kx_j|=0$ for all
$m>k$, consequently, $1\ne \sum_{j=1}^kx_j$ for each $x\in A_1^*$
and $k\in {\bf N}$, since $|1-b|=0$ is equivalent to $b=1$ and $x$
does not belong to $\bigcup_nE^n$. \par We take an arbitrary vector
$x=(x_1,x_2,...)$ in $A_1^*$. Since $x_1\in B({\bf K},0,1)\setminus
\{ 1 \} $, we can put $y_1=x_1$. Moreover, we have
\par $(1)$ $|1-x_1-...-x_k|\le \max (|1-x_1 -...- x_{k+1}|, |x_{k+1}|)$ and
\par $(2)$ $|x_{k+1}| \le \max (|1-x_1-...-x_{k+1}|, |1-x_1-...-x_k|)=
|1-x_1-...-x_k|$ \\ for each $k$ and each $x\in A_1$. So we can take
$y_2=1-(1-x_1-x_2)/(1-x_1)=x_2/(1-x_1)$. By induction if
$x_1,...,x_m$ are marked, then $x_{m+1}\in B({\bf
K},-x_1-...-x_m,1)\setminus \{ 1 \} $ and it is sufficient to put
\par $(3)$ $y_{m+1} =
1-(1-x_1-...-x_{m+1})/(1-x_1-...-x_m)=x_{m+1}/(1-x_1-...-x_m)$,
consequently, $y_j\in B({\bf K},0,1)$ for each $j\in {\bf N}$.
Therefore,
\par $(4)$ $(1-x_1-...-x_{m+1})=(1-y_{m+1})(1-x_1-...-x_m)=(1-y_{m+1})...(1-y_1)$
for each natural number $m\in {\bf N} = \{ 1,2,3,... \} $. The
inverse mapping $q^{-1}$ is given by: \par $(5)$ $x_1=y_1$,
$x_2=y_2(1-y_1)$, $x_{m+1} = y_{m+1}(1-x_1-...-x_m) =
y_{m+1}(1-y_m)...(1-y_1)$ for every natural number $m\in {\bf N}$,
consequently, $|x_j|\le 1$ for each $j$ and $|\sum_{j=1}^kx_j|\le
\max_{j=1}^k|x_j|\le 1$ for each $k\in {\bf N}$. Thus $x_{m+1}=0$ is
equivalent to $y_{m+1}=0$, since $y_j\ne 1$ for each $y\in s$. But
for each $y\in s_*$ the set $\{ j: ~ y_j \ne 0 \} $ is infinite,
which is equivalent to the fact that the set $ \{ j: ~ x_j\ne 0 \} $
is infinite for $x=q^{-1}(y)$. That is, $q^{-1}(s_*)\subset A_1^*$.
\par From the definition of $q$ one can lightly see
that $\lim_{k\to \infty }\sum_{j=1}^kx_j=1$ is equivalent to the
fact that the limit \par $(6)$ $\lim_{j\to \infty }
(1-y_j)...(1-y_1)=0\in B({\bf K},0,1)$ exists. Thus the mapping $q$
is bijective from $A_1^*$ onto $s_*$, moreover, $q$ and its inverse
$g=q^{-1}: s_*\to A_1^*$ are coordinate-wise continuous.
\par Let $x^n$ be a converging sequence in $A_1^*$, $\lim_{n\to
\infty } x^n =x\in A_1^*$, then \par $y_j^{n+m} - y_j^n =
x_j^{n+m}/(1-x_1^{n+m}-...-x_{j-1}^{n+m}) -
x_j^n/(1-x_1^n-...-x_{j-1}^n)$ and \par $|y_j^{n+m}-y_j^n| \le $
\par $\max_{0\le i \le j-1}
|x_j^{n+m}x_i^n-x_j^nx_i^{n+m}|/[|1-x_1^{n+m}-...-x_{j-1}^{n+m}|
|1-x_1^n-...-x_{j-1}^n|]\le $ \\
$\max_{0\le i \le j-1} (|x_j^{n+m}-x_j^n| |x_i|^n, |x_j^n|
|x_i^{n+m}-x_i^n|)/[|1-x_1^{n+m}-...-x_{j-1}^{n+m}|
|1-x_1^n-...-x_{j-1}^n|],$ \\ where $x_0^n=1$, consequently, the
mapping $q: A_1^*\to s_*$ is continuous, since $s_*$ is in the
topology inherited from the Tychonoff product topology on $s$ (see
\S 7).
\par If $y^n$ is a converging sequence in $s_*$, then
\par $x_j^{n+m}-x_j^n= y_j^{n+m}(y_{j-1}^{n+m}-1)...(y_1^{n+m}-1)-
y_j^n(y_{j-1}^n-1)...(y_1^n-1)$, consequently,
\par $(7)$ $|x_j^{n+m}-x_j^n|\le \max ( |y_j^{n+m}-y_j^n|
|(y_{j-1}^{n+m}-1)...(y_1^{n+m}-1)|, $\par $ ~ |y_j^n|
|(y_{j-1}^{n+m}-1)...(y_1^{n+m}-1)  - (y_{j-1}^n-1)...(y_1^n-1)|)$
and \par $\| x^{n+m} - x^n \| = \sup_j |x_j^{n+m}-x_j^n| $. \\ But
$|x_j|\le 1$ for each $j$ and $\lim_{j\to \infty }x_j=0$. For each
$\epsilon
>0$ there exists a natural number $j_0>0$ such that $|x_j|<\epsilon $ for each
$j>j_0$. Then for each $\delta
>0$ there exists a natural  number $n_0$ such that
$|y_k^n-y_k|<\delta $ for each $k=1,...,j_0$. Choose $\delta >0$
such that $|x_k-x_k^n|<\epsilon $ for each $k=1,...,j_0$ and
$n>n_0$. Then \par $|x_j|\le \max (|1-x_1-...-x_j|,
|1-x_1-...-x_{j-1}|)\le |1-x_1-...-x_{j_0}|$
\par $ =|(1-y_{j_0})...(1-y_1)|$ and \par $(8)$ $|x_j^n|\le
|(1-y_{j_0}^n)...(1-y_1^n)|$ \\ for all $n\in {\bf N}$ and each
$j>j_0$.  \\ Therefore, from Formulas $(7,8)$ it follows that the
mapping $g=q^{-1}: s_*\to A_1^*$ is continuous, since $(1-y_k)\in
B({\bf K},0,1)$ for each $k$, $|ab|\le |a| |b|$ for all $a, b\in
{\bf K}$, while \par $|x_j-x_j^n|\le \max (|x_j|,|x_j^n|)\le
|(1-y_{j_0})...(1-y_1)|(1+\delta )$ \\ for each $n>n_0$ and $j>j_0$,
where $\delta >0$ can be taken less than $\epsilon $.
\par {\bf 11. Lemma.} {\it Each element $x\in A_2$ can be presented
as $x=x^1+x^2$, where $1-x^1$ and $x^2\in A_1$.}
\par {\bf Proof.} Consider the bounded sequence
\par $a_n := |1-\sum_{j=1}^kx_j|\le \max (1,|\sum_{j=1}^kx_j|)=1$ \\
in ${\cal R}$ and compose the new sequence $w_1=a_1$, $w_n =a_n-
a_{n-1}$ for each $2\le n\in {\bf N}$, consequently,
$a_n=w_1+...+w_n$ for each $2\le n\in {\bf N}$. Then we put $b_n :=
\max (w_n,0)$ and $c_n := - \min (w_n,0)$, consequently,
$w_m=b_m-c_m$ and $0\le b_m$ and $0\le c_m$ for each $m\in {\bf N}$.
Since $\sum_{j=1}^kx_j\ne 1 ~ \forall k$, and $~ \sum_{j=1}^{\infty
}x_j=1$, we certainly have the conditions $a_n> 0$ for each $n\in
{\bf N}$ and $\lim_n a_n=0$.
\par Put now $d_n := \sup_{m\ge n} [\max (b_m,c_m)+ a_{m-1}]$ and $e_n :=
\inf_{1\le m\le n} a_m$, where $a_0:=0$, consequently, $0\le
d_{n+1}\le d_n\le 2$ and $0\le e_{n+1}\le e_n\le 1 $ for each $n\in
{\bf N}$ and $\lim_n d_n=0$ and $\lim_n e_n=0$, since
$$||1-\sum_1^{m-1} x_j| - |1-\sum_1^m x_j|| \le |x_m|\le 1$$ for each
$m\in {\bf N}$ and $x\in A_2$. From $a_n \le \max (b_n+a_{n-1}, ~
c_n + a_{n-1})$ one gets the inequality $a_n\le d_n$ for each
natural number $n$. On the other hand, $e_n\le a_n$ and $a_n, e_n
\in \Gamma _{\bf K}\bigcup \{ 0 \} $ for each $n$.
\par If $y, z \in {\bf K}$, $|y| =r_1\le r_2$, $|z|=r_2$, then
$(r_2-r_1) \le |y+z| \le r_2$.
\par Consider subsets $ \{ (\beta _n,\gamma _n)\in {\bf K}^2: ~ |\beta
_n | =a_n, ~ |1-\gamma _n| = e_n, ~ |1-\beta _n - \gamma _n |=a_n \}
$. Therefore, $1-\beta _n =x^1_1+...+x^1_n$ and $\gamma _n =
x^2_1+...+x^2_n$ give two elements $1-x^1$ and $x^2\in A_1$ such
that $x=x^1+x^2$, where $x^l=(x^l_1,x^l_2,...)$, $x^l_k\in {\bf K}$
for each $k\in {\bf N}$ and $l=1, 2$.
\par {\bf 11.1. Corollary.} {\it There are embeddings:
\par $(1)$ $A_2\hookrightarrow (1-A_1)\cup A_1\hookrightarrow c_0$
and
\par $(2)$ $A_2^*\hookrightarrow (1-A_1^*)\oplus A_1^*\hookrightarrow
c_0$.}
\par {\bf Proof.} The first embedding follows from Lemma 10.
On the other hand, $(1-A_1^*)\cap A_1^*=\emptyset $, since
$\sum_{j=1}^k x_j\ne 1$ for each $x\in A_1^*$ and $k\in {\bf N}$
while $\sum_{j=1}^{\infty }x_j=1$, but $\sum_{j=1}^{\infty } y_j=0$
for each $y\in 1-A_1^*$. Therefore, $(1-A_1^*)\cup
A_1^*=(1-A_1^*)\oplus A_1^*$, since the mapping $x\mapsto
\sum_{j=1}^{\infty } x_j$ is continuous from $c_0$ into $\bf K$, at
the same time $$|\sum_{j=1}^{\infty } x_j -\sum_{j=1}^{\infty } z_j
|\le \sup_{j\in {\bf N}} |x_j-z_j| = \| x - z \| $$ for every $x, z
\in c_0$.
\par {\bf 12. Lemma.} {\it The topological spaces $A_1^*$ and
$A_2^*$ are homeomorphic.}
\par {\bf Proof.} The Banach spaces $c_0$ and $c$ are linearly topologically isomorphic,
since $c_0\oplus {\bf K}$ is linearly topologically isomorphic with
$c_0$ and with $c$ as well, where $c=c(\omega _0,{\bf K})$,
$$(1)\quad c(\alpha ,{\bf K}) := \{ x: ~ x= (x_j: j\in \alpha ), \forall
j\in \alpha ~ ~ x_j\in {\bf K},$$
$$\| x \| =\sup_{j\in \alpha } |x_j|, \exists \lim_j x_j\in {\bf K}
\}.$$ If $y\in c_0$ put \par $(2)$ $x_1=y_1$, $x_j=y_1+...+y_j$ for
each $2\le j\in \omega _0$. \par Since $\lim_j y_j=0$ and
$|\sum_{j=n}^m y_j|\le \max_{n\le j\le m} |y_j|$ for each $1\le n
\le m$, the series $\sum_{j=1}^{\infty } y_j$ converges and $x\in
c$. Consider $A_1^*$ and $A_2^*$ embedded into the Banach space $c$,
$A_l^*\ni y\mapsto x=x(y)\in c$ (see Formula $(2)$). Take $x\in
A_l^*$, then $ \| x \| = \epsilon >0$, where $l=1,2$. There exists
$\xi \in A_l^*$ with $ \| x \| = \| \xi \| $ and $\inf_j |\xi _j
|=\delta >0$ choosing $0<\delta \le \min (1/p,\epsilon )$.
Therefore, if $\xi \in A_1^*$, then $B(A_2^*,\xi ,\delta /p) \subset
A_1^*$. Indeed, if $|a-x_j|\le \delta /p$ and $|b-x_{j+1}| \le
\delta /p$, then $|a| = |x_j|$ and $|b|=|x_{j+1}|$, where $a, b \in
{\bf K}$. That is from $|x_{j+1}|\le |x_j|$ it follows, that $|b|\le
|a|$. \par Let $\omega ({\bf K})$ be the topological weight of the
field $\bf K$, then $A_1^*$ and $A_2^*$ have coverings by disjoint
families ${\cal F}_1$ and ${\cal F}_2$ of balls such as $B(A_2^*,\xi
,r)$ with $0<r\le \delta /p$ as above, where $| {\cal F}_1 | = |
{\cal F}_2 |= \omega ({\bf K}) |\Gamma _{\bf K}| \aleph _0,$ as
usually $|E|=card (E)$ denotes the cardinality of a set $E$. Each
two balls $B(A_2^*,\xi ,r_1)$ and $B(A_2^*,\eta ,r_2)$ in $A_2^*$ of
radius less than $1/p$ are homeomorphic, since they are disjoint
unions of $\omega ({\bf K})\aleph _0$ balls $B(A_2^*,z,r_3)$ of
radius $r_3=\min (r_1,r_2)/p$. Hence $A_1^*$ and $A_2^*$ are
homeomorphic due to Corollary 11.1.
\par {\bf 13. Definitions.} A topological space $X$ is called
ultrametric, if its topology is given by an ultrametric $\rho $
having values in $\Gamma _{\bf K}\cup \{ 0 \} $ such that
\par $(1)$ $\rho (x,y)\ge 0$ for every $x, y\in X$ and
$\rho (x,y)=0$ if and only if $x=y$;
\par $(2)$ $\rho (x,y)=\rho (y,x)$ for every $x, y \in X$;
and a metric satisfies the ultrametric inequality which is stronger
than the usual triangle inequality:
\par $(3)$ $\rho (x,z) \le \max (\rho (x,y), \rho (y,z))$ for every
$x, y, z\in X$.
\par If $X$ is an ultrametric space and $H_t: X\times B({\bf
K},0,1)\to X$ is a simultaneously continuous mapping so that $t\in
B({\bf K},0,1)$ and $H_0=id$, $~H_t: X\to X$ is a homeomorphism for
each $t$, then $H_t$ is called an isotopy. The isotopy is called
invertible, if $H_t^{-1}$ is jointly continuous in $(x,t)\in X\times
B({\bf K},0,1)$.
\par {\bf 14. Lemma.} {\it Let $(X,\rho )$ be an ultrametric space, and let $A$ and
$B$ be two non intersecting closed subsets in $X$ such that
$\inf_{x\in A, ~ y \in B} \rho (x,y)= d >0$. Suppose that $\bf K$ is
an infinite non discrete non-archimedean field. Then there exists a
continuous function $f: X\to {\bf K}$ such that $f(A)=\{ 0 \} $ and
$f(B)= \{ 1 \} $.}
\par {\bf Proof.} Since $A$ and $B$ are closed subsets in $X$
and $A\cap B = \emptyset $, then $\inf_{x\in A, ~ y \in B} \rho
(x,y)= d >0$. A field $\bf K$ is non discrete, that is $0$ is a
limit point of a set $ \{ |x|: x\ne 0, ~ x\in {\bf K} \} $ in a ring
${\cal R}$. Therefore, an element $r\in \Gamma _{\bf K}$ exists such
that $0<r< d$. We take two subsets $U:= \{ x\in X: ~ \rho (B,x) \le
r \} $ and $V := \{ x\in X: ~ \rho (A,x) \le r \} $, where $\rho
(B,x) := \inf_{y\in B} \rho (x,y)$. From the ultrametric inequality
it follows that two subsets $U$ and $V$ do not intersect, $U\cap
V=\emptyset $. The set $U$ is clopen in $X$, since $B(X,y,r)\subset
U$ for each $y\in U$, where $B(X,y,r) := \{ z\in X: ~ \rho (y,z)\le
r \} $, consequently, $X\setminus U$ is clopen in $X$. On the other
hand, $A\subset V\subset (X\setminus U)$. Next we take a function
$f: X\to {\bf K}$ such that $f(x)=1$ for each $x\in U$ and $f(x)=0$
for each $x\in X\setminus U$. This function $f$ is continuous
\par {\bf 15. Lemma.} {\it Let $A_1,...,A_n$ be closed non intersecting
(i.e. disjoint) subsets in an ultrametric space $(X,\rho )$ such
that $\min_{k\ne l} \inf_{x\in A_k, ~ y \in A_l} \rho (x,y)= d >0$.
Suppose that $\bf K$ is an infinite non discrete non-archimedean
field and $b_k\ne b_j\in {\bf K}$ for each $k\ne j$, $j, k=1,...,n$.
Then there exists a continuous function $f: X\to {\bf K}$ such that
$f(A_k) = \{ b_k \} $.}
\par {\bf Proof.} We choose $r\in \Gamma _{\bf K}$ such that
$0<r<d$ and take clopen subsets $U_j = \{ x\in X: ~ \rho (A_j,x)\le
r \} $ for each $j=1,...,n$. Then $A_j\subset U_j$ for each $j$ and
$A_n\subset X\setminus (\bigcup_{k=1}^{n-1} U_j)$. These clopen
subsets $U_j$ are pairwise disjoint due to the ultrametric
inequality. The field ${\bf K}$ is infinite, hence there are
distinct elements $b_1,...,b_n$, that is $b_k\ne b_j\in {\bf K}$ for
each $k\ne j$, $j, k=1,...,n$. Take any distinct $n$ elements
$b_1,...,b_n$ in ${\bf K}$. We construct a function $f: X\to {\bf
K}$ such that $f(x)=b_j$ for each $x\in U_j$ with $j=1,....,n-1$ and
$f(x)=b_n$ for each $x\in X\setminus (\bigcup_{k=1}^{n-1} U_j)$.
Since subsets $U_j$ for each $j$ and $X\setminus
(\bigcup_{k=1}^{n-1} U_j)$ are clopen in $X$, then the function $f$
is continuous.

\par {\bf 15.1. Lemma.} {\it Let $A$ and $B$ be two balls in an ultrametric space
$(X,\rho )$ such that $\inf_{x\in A, ~ y \in B} \rho (x,y)= 0$. Then
one of the balls is contained in another.}
\par {\bf Proof.} Let $r_a$ and $r_b$ be radiuses of balls $A$ and
$B$ such that $r_b\le r_a$. Since $\inf_{x\in A, ~ y \in B} \rho
(x,y)= 0$, then for each $\epsilon >0$ in ${\cal R}$ there exist
$x\in A$ and $y\in B$ such that $\rho (x,y)<\epsilon $. For every
$0<\epsilon <r_b$, $z\in B$ and $v\in A$ this gives $$\rho (x,z)\le
\max (\rho (x,y), \rho (y,z))\le r_b \mbox{ and}$$
$$\rho (v,z)\le \max (\rho (v,x), \rho
(x,z))\le r_a,$$ consequently, $z\in A$ and hence $B\subset A$.

\par {\bf 15.2. Lemma.} {\it Suppose that $(X,\rho )$ is an
ultrametric space, $A$ and $B$ are two disjoint subsets in $X$,
$A\cap B=\emptyset $, such that $A$ and $B$ are unions of balls of
radiuses $r$ with $r\ge d$, where $d\in {\cal R}$, $d>0$. Then
$\inf_{x\in A, y\in B} \rho (x,y)\ge d$.}
\par {\bf Proof.} If $\inf_{x\in A, y\in B} \rho (x,y)< d$, then two
pints $x_0\in A$ and $y_0\in B$ exist such that $\rho (x_0,y_0)<d$,
since otherwise $\rho (x,y)\ge d$ for all $x\in A$ and $y\in B$.
Then we take balls $B_1$ in $A$ and $B_2$ in $B$. Their radiuses
$r_1$ and $r_2$ are not less than $d$, $r_1\ge d$ and $r_2\ge d$.
Therefore, for each $x\in B_1$ and $y\in B_2$ the inequalities are
satisfied:
$$\rho (x,y) \le \max (\rho (x,x_0), \rho (x_0,y_0), \rho (y_0,y))\le
\max (r_1,r_2),$$ consequently, either $B_1\subset B_2$ for $r_1\le
r_2$ or $B_2\subset B_1$ for $r_2\le r_1$. Thus $A\cap B\ne
\emptyset $ contradicting the suppositions of this lemma.

\par {\bf 16. Theorem.} {\it  Let $(X,\rho )$ be an ultrametric space.
Suppose that ${\bf K}$ is a non-archimedean infinite non discrete
locally compact field. Let also $M$ be a closed subset in $X$ and
let $f: M\to {\bf K}$ be a continuous function. Then there exists a
continuous extension $g$ of $f$ on $X$, $~g: X\to {\bf K}$.}
\par {\bf Proof.} In view of Lemma 6 ${\bf K}$ and $B({\bf K},0,1,)\setminus \{ 1 \} $
are homeomorphic. Therefore, it is sufficient to consider continuous
mappings $f: M\to B({\bf K},0,1)$. In this case a constant $c\in
\Gamma _{\bf K}$ exists such that $|f(x)|\le c$ for each $x\in M$.
On the other hand, each locally compact field ${\bf K}$ has a
discrete group $\Gamma _{\bf K}$.
\par We construct a sequence of functions $g_n: X\to {\bf K}$ a
limit of which $g: X\to {\bf K}$ will be a continuous extension of
$f$. Consider a continuous function $g_n|_M: M\to {\bf K}$ from a
closed subset $M$ in $X$. A field ${\bf K}$ is locally compact,
consequently, the quotient ring $B({\bf K},0,c)/B({\bf K},0,c/p^n)$
is finite. Hence there are points $x_{n,j}\in B({\bf K},0,c)$ such
that $B({\bf K},0,c)$ is a disjoint union of clopen balls $B({\bf
K},x_{n,j},c/p^n)$, where $ j=1,...,m(n)\in {\bf N}$. We take the
clopen subsets $A_{n,j} := (g_n|_M)^{-1}(B({\bf K},x_{n,j},c/p^n))$
in $M$ and consider clopen subsets $$V_{n,j} := \{ x\in X: ~ \rho
(x,A_{n,j})\le c/p^{n+1} \} $$ in $X$ and then put
$$U_{n,1}=V_{n,1}\mbox{ and  }U_{n,j} = V_{n,j}\setminus \bigcup_{l<j} V_{n,l}$$ for each $j>1$.
Thus $\bigcup_{j=1}^{m(n)} A_{n,j} = M$. Since $g_n$ is continuous,
then from $\lim_k \rho (x_k,y)=0$ in $X$ it follows that $\lim_k \|
g_n(x_k)-g_n(y) \| =0$.
\par Each clopen subset $U_{n,j}$ is the disjoint union of balls in
$X$. Let $W_{n,j,k}$ be the disjoint union of balls in $U_{n,j}$ of
radiuses $r$ with $r\ge p^{-k}$, where $k\in {\bf N}$, hence
$W_{n,j,k}\subset U_{n,j}$ for each $k\in {\bf N}$. Then by Lemma
15.2 the inequality is valid:
$$\inf_{j\ne l} \inf_{x\in W_{n,j,k}, y\in W_{n,l,k}} \rho (x,y)\ge
p^{-k}.$$ If $\{ f_{n,k}: k \in {\bf N} \} $ is a family of
continuous functions on $\bigcup_jW_{n,j,k}$ such that their
restrictions satisfy the condition $$f_{n,k+1}|_{\bigcup_jW_{n,j,k}}
= f_{n,k}|_{\bigcup_jW_{n,j,k}}$$ for each $k\in {\bf N}$, then
there exists their continuous combination $f_n:=\nabla_{k\in {\bf
N}} f_{n,k}$ on $\bigcup_{k=1}^{\infty } \bigcup_jW_{n,j,k} =
\bigcup_j U_{n,j}$, since $\bigcup_jW_{n,j,{k+1}}\setminus
(\bigcup_jW_{n,j,k})$ is the disjoint union of balls of radiuses $r$
with $p^{-k-1}\le r<p^{-k}$.
\par In view of Lemmas 15, 15.1 and the proof above a continuous
function $u_n: X\to {\bf K}$ exists satisfying the conditions
\par $(1)$ $u_n(U_{n,j})\subset B({\bf K},x_{n,j},c/p^n)$ for each
$j=1,...,m(n)$ and \par $(2)$ $|u_n(x)|\le c$ for each point $x$ in
$X$. Therefore, \par $(3)$ $|u_n(x)-g_n(x)|\le c/p^n$ for each $x\in
M$ and \par $(4)$ $M\subset \bigcup_{j=1}^{m(n)} U_{n,j}=: P_{n+1}$.
\par Each subset $P_n$ is clopen in $X$ and it is possible to take any
continuous function $h: X\setminus P_2\to {\bf K}$. We continue this
process by induction from $n=1$ with $g_1|_M=f$, further putting
\par $(5)$ $g_{n+1}|_M=f: M\to {\bf K}$ and \par $(6)$ $g_{n+1}(y)=u_n(y)$ for each
$y\in X\setminus (P_n\setminus P_{n+1})$ when $n\ge 2$ and
\par $(7)$ $g_{n+1}(y)=h(y)$ for each $y\in X\setminus P_2$. \par Then
$P_{n+2}\subset P_{n+1}$ and $P_{n+2}$ is clopen in $P_{n+1}$ for
each $n\in {\bf N} = \{ 1, 2, 3,... \} $ and
\par $(8)$ $\bigcap_{n=2}^{\infty } P_n=M$. \par Therefore, $g_{n+k}=g_{n+1}$ on
$X\setminus P_{n+2}$ for each $k\ge 2$. From the construction above
and formulas $(1-8)$ it follows, that $\lim _n g_n(x)=f(x)$ for each
$x\in M$ and $g=\lim_n g_n$ is continuous from $X$ into $B({\bf
K},0,c)$, since $|g_n(x) - g_{n+1}(x)| \le cp^{-n}$ for each $x\in
X$ and $n\in \bf N$.
\par {\bf 17. Example.} Let $f_j $ be a continuous surjective
mapping  from $B({\bf K},0,1)$ into $B({\bf K},0,1)\setminus B({\bf
K},0,1/p^j)$ such that $f_j(B({\bf K},0,1)\setminus B({\bf
K},0,1/p^j)) = \{ 1 \} $ and $f_j(B({\bf K},0,1/p^{j+1})) = \{ 1 \}
$ and $f_j(B({\bf K},0,1/p^j)\setminus B({\bf K},0,1/p^{j+1}))= \{
b_j \} $, $ ~ b_j\ne 0$ with $\lim_j b_j=0$. Next a mapping $H_t (x)
=(x_1f_1(t), x_2f_2(t),...)$ exists for each $x=(x_1,x_2,...,)\in
c_0$ and $t\in B({\bf K},0,1)$. Then $H_t$ is an isotopy, but
$H_t^{-1}$ is not jointly continuous at $(0,0)\in c_0\times B({\bf
K},0,1)$, since $B({\bf K},0,1)=B({\bf K},1,1)$. On the other hand,
$H_0^{-1}$ is continuous on $c_0$ and $H_t^{-1}(0)$ is continuous in
$t\in B({\bf K},0,1)$.
\par {\bf 18. Definition.} For a subset $K$ of a totally
disconnected Hausdorff topological space $X$ let $H_t$ be a
$1$-parameter family of homeomorphisms onto $X$, with $t\in B({\bf
K},0,1)$, so that $H_0=id$, $~H_t(X)=X$ for each $t\ne 1$,
$H_1(X\setminus K) =X$ and $H_t$ and $H_t^{-1}$ are jointly
continuous in $(x,t)\in X\times B({\bf K},0,1)$. Then $H_t$ is
called an invertible non-archimedean isotopy pushing $K$ off $X$.
Particularly, $K$ may be a singleton $ \{ p \} $.
\par If $a\ne b\in B({\bf K},0,1)$ with $0\le |a| \le 1$ and $0\le
|1-b| <1$, we put $H[a,b]_t=id$ for each $t\in B({\bf
K},0,|a|)\setminus B({\bf K},1,|1-b|)$, $~H[a,b]_t = H_1$ for $t\in
B({\bf K},1,|1-b|)$, $~H[a,b]_t = H_{(t-a)/(b-a)}$ for each $t\in
B({\bf K},0,1)\setminus [B({\bf K},0,|a|)\cup B({\bf K},1,|1-b|)]$.
If $a,b ,c \in B({\bf K},0,1)$, $~ 0\le |a|\le 1$, $~0<|1-b|<1$,
$~0\le |1-c|\le |1-b|$, $~a\ne b$, $~b\ne c$, the composition
$F[b,c]_t\circ H[a,b]_t$ of two isotopies is defined.
\par {\bf 19. Lemma.} {\it Let $H^j$ be an invertible (non-archimedean) isotopy
of a totally disconnected Hausdorff topological space $X$ onto $X$
for each $j\in {\bf N}$ and let $p_0\in X$ be a marked point. Let
also
\par $(1)$ for each point $x\in X\setminus \{ p_0 \} $ a
neighborhood $U$ of $x$ in $X$ and a natural number $n(U)\in {\bf
N}$ exist so that $H_t^n=id$ on $H_1^{n(U)}\circ ... \circ
H^2_1\circ H^1_1(U)$ for each $n>n(U)$ and \par $(2)$ for each point
$y\in X$ a neighborhood $V$ of $y$ in $X$ and an integer $n(V)$
exist so that $(H_t^n)^{-1}=id$ on $V$ and $H_1^{n(V)}\circ ...\circ
H_1^2\circ H_1^1(V)\subset (X\setminus \{ p_0 \} )$ for each
$n>n(V)$.
\par Then $H_t=L\prod_{j=0}^{\infty } H^{j+1}[1-p^j,1-p^{j+1}]_t$ is an
invertible (non-archimedean) isotopy pushing $p_0$ off $X$.}
\par {\bf Proof.} For each $t\in B({\bf K},0,1)$, $~|t|\ne 1$, a natural number
$k\in \bf N$ exists such that $t\in B({\bf K},0,p^{-k})$,
consequently, $H^j[1-p^j,1-p^{j+1}]_t=id$ on $X$, when $j\ge 1$,
since $|t|<|1-p^j|=1$. If $|t|=1$ and $t\ne 1$ a non-negative
integer $k$ exists such that $t\in B({\bf K},1,p^{-k})\setminus
B({\bf K},1,p^{-k-1})$. Therefore, $H_t$ reduces to a finite product
of homeomorphisms of $X$ onto $X$. This implies that $H_t$ and
$(H_t)^{-1}$ are continuous for each $t\in B({\bf K},0,1)$, $ ~ t\ne
1$. \par If $y\in X\setminus \{ p_0 \} $ and $t=1$, the continuity
of $H$ at $(y,1)$ follows from Condition 1. On the other hand,
Condition 2 implies the continuity of $H^{-1}$ at $(y,1)$ for any
$y\in X$.
\par {\bf 20. Lemma.} {\it Consider $\bf s$ as in \S 3. Let either $X=\bf s$
or $X=c_0$. Suppose that $T$ is a compact subset in $X$. Then a
homeomorphism $\eta $ of $X$ onto $X$ exists so that $\pi _1\circ
\eta (T)$ is a single point in ${\bf K}$, where $\pi _j: X\to {\bf
K}_j$ is a projection linear over $\bf K$.}
\par {\bf Proof.} The topological vector spaces $\bf s$ and $c_0$
have projections $\pi _j$ linear over $\bf K$. The demonstration
below is given for $\bf s$, whilst that of $c_0$ is analogous. Take
homeomorphisms $f_j$ of ${\bf K}_1\times {\bf K}_j$ onto itself,
where ${\bf K}_j={\bf K}$ for each $j\in {\bf N}$, satisfying two
conditions:
\par $(1)$ $f_j(x_1,x_j) = (x_1,y_j)$ for each $(x_1,x_j)\in {\bf
K}_1\times {\bf K}_j$, where $y_j\in {\bf K}_j$, and
\par $(2)$ if $D_j$ is a region in ${\bf K}_1\times {\bf K}_j$
such that $D_j = \{ (x,y): ~ (x,y)\in {\bf K}_1\times {\bf K}_j, ~
|x-x_0| \le a, ~ |y-y_0|\le b \} $ for some $(x_0,y_0)\in {\bf
K}_1\times {\bf K}_j$, $~a, b \in \Gamma _{\bf K}$, $~ \pi
_{1,j}(T)\subset D_j$, where $\pi _{1,j}: {\bf s}\to {\bf K}_1\times
{\bf K}_j$ is a linear over ${\bf K}$ projection, and $f_j(D_j) = \{
v: ~ v\in {\bf K}_1\times {\bf K}_j, ~ v_1=
t_1(b_1-a_1)+t_2(c_1-a_1), ~  v_2 = t_1(b_2-a_2) + t_2(c_2-a_2), ~
t_1, t_2\in B({\bf K},0,1) \} $, $~ a=(a_1,a_2)$, $~b=(b_1,b_2)$, $~
c= (c_1,c_2)$ are marked points in ${\bf K}_1\times {\bf K}_j$ and
$v=(v_1,v_2)$ such that $f_j(D_j)\cap (w+({\bf K},0)) = E_{w,j}$
with \par $diam E_{w,j} := \sup_{x, y \in E_{w,j}} |x-y|\le p^{-j}$
\\ for each $w \in {\bf K}_1\times {\bf K}_j$.
\par Then a homeomorphism of $f$ of $\bf s$ onto itself exists such
that $f(x_1,x_2,...) = (x_1,y_2,y_3,...)$, where $y_j$ is given by
Condition $(1)$ for each $j$. From the construction of $f$ it
follows that $f$ is bijective and surjective. Since each $f_j$ is
continuous  and $\bf s$ is supplied with the Tychonoff topology,
then the mapping $f$ is also continuous.
\par From Condition $2$ it follows that $f(T)\cap (v+({\bf K},0))$
is  either a singleton or the void set for each $v\in {\bf s}$.
\par Put ${\bf s}_0 := \{ x: ~ x \in {\bf s}, ~ x_1=0 \} $ and let
$\pi _0: {\bf s}\to {\bf s}_0$ be the corresponding linear over
${\bf K}$ projection onto ${\bf s}_0$. Therefore, in accordance with
the construction above $\pi _0|_{f(T)}$ is a homeomorphism from
$f(T)$ into ${\bf s}_0$. Consider the restriction $\Phi = \pi
_1\circ \pi _0^{-1}|_{\pi _0(f(T))}$. In view of Theorem 16 it has a
continuous extension $\psi : {\bf s}_0\to {\bf K}$, that is, $\psi
|_{\pi _0(f(T))} = \Phi $.
\par There exists a homeomorphism $\xi $ of ${\bf s}$ onto ${\bf s}$
so that $\xi |_{(p+({\bf K},0))} (x) = x - \Phi (p)$ for each $p\in
{\bf s}_0$ and $x\in p+({\bf K},0)\subset \bf s$. The desired
homeomorphism is $\eta = \xi \circ f$, since $\eta (T)\subset {\bf
s}_0$ and $\pi _1(\eta (T)) = \{ 0 \} $.
\par {\bf 21. Theorem.} {\it  Let $ \{ T_j: ~ j\in {\bf N} \} $ be a family
of compact subsets of $X={\bf s}$ or $X=c_0$. Then a homeomorphism
of $X$ onto $X$ exists such that each subset $g(T_j)$ is infinitely
deficient for each $j\in {\bf N}$.}
\par {\bf Proof.} If $\beta _j$ is a family of disjoint subsets of
$\bf N$ so that $\bigcup_j \beta _j=\bf N$, then $\bf s$ and $c_0$
can be written as ${\bf s}=\prod_j {\bf s}^{\beta _j}$ (see \S 3),
and $c_0=c_0(c_0(\beta _j): j)$ respectively, where \par $c_0(Y_j:
j) = \{ y=(y_1,y_2,...): ~ \forall j ~ y_j\in Y_j, \forall \epsilon
>0$\par $ card \{ j: ~ \| y_j \| >\epsilon \} < \aleph _0; ~ \| y \|
=\sup_j \| y_j\|_{Y_j} \} $ \\ for a family of Banach spaces $Y_j$
over $\bf K$. Suppose that a mapping $\theta : {\bf N}\to {\bf N}$
has the property that $\theta ^{-1} (j)$ is infinite for each $j\in
{\bf N}$. Then each ${\bf s}^{\beta _j}$ or $c_0(\beta _j)$
respectively with $\beta _j = \theta ^{-1}(j)$ is homeomorphic with
$\bf s$ or $c_0$ respectively. In view of Lemma 20 a homeomorphism
$g_j$ of ${\bf s}^{\beta _j}$ or $c_0(\beta _j)$ onto itself so that
$g_j(\pi _{\beta _j}(T_{\theta (j)}))$ is deficient relative to the
first element of $\beta _j$. Put $g=\prod_j g_j$ for $\bf s$ or
$g(y)=(g_1(y_1),g_2(y_2),...)$ for each $y\in c_0=c_0(c_0(\beta _j):
j)$ respectively. Then $g(T_j)$ is infinitely deficient for each
$j\in \bf N$.
\par {\bf 22. Remark.} Let $\mbox{}_rH_t$ be a two-parameter
(non-archimedean) family of homeomorphisms with $r\in B({\bf
K},0,1)\setminus \{ 0 \} $ and $t\in B({\bf K},0,1)$ such that for
$r$ fixed $\mbox{}_rH$ is an isotopy pushing the origin off $\bf s$.
If $\mbox{}_rH_tx$ and $(\mbox{}_rH_t)^{-1}x$ are continuous in $r,
t$ and $x\in \bf s$ or $c_0$, then the family $ \{ \mbox{}_rH_t: ~ r
\in B({\bf K},0,1)\setminus \{ 0 \} \mbox{ and } t\in B({\bf K},0,1)
\} $ is called an invertible continuous family of invertible
(non-archimedean) isotopies. Henceforth, the case is considered,
when $\mbox{}_rH_t$ is the identity outside the $|r|$-neighborhood
of the origin in $\bf s$ or $c_0$ respectively. The topological
space ${\bf s}$ is metrizable with the complete metric
\par $(1)$ $d(x,y) = \sum_j \min (p^{-j}, |x_j-y_j|^{-1})\in {\cal R}$.
\par {\bf 23. Lemma.} {\it There exists an invertible
(non-archimedean) isotopy $H$ pushing the origin $x_0$ off $X$,
where $X=\bf s$.}
\par {\bf Proof.} We consider an invertible (non-archimedean)
isotopy $F^j$ on ${\bf K}_1\times {\bf K}_{j+1}$ satisfying the
following conditions:
\par $(1)$ $F_t^j (x_1,x_{j+1}) =
(x_1,x_{j+1})$ for each $(x_1,x_{j+1})\in {\bf K}_1\times {\bf
K}_{j+1}$ with $|x_1|<p^{j-1}$,
\par $(2)$ $F_1^j$ maps the set $A_j := \{ (z_1,z_2): ~ z_1=\xi _1,
z_2=\xi _2\gamma +\xi _3(1-\gamma ), \gamma \in B({\bf K},0,1) \} $
with $|\xi _1| = p^{j-1}$ and $|\xi _2|=p^{-j}$ and $|\xi _3|=p^j$
onto the set $B_j := \{ (z_1,0): ~ z_1=(\xi _2+\xi _4)\gamma + (\xi
_3+\xi _4)(1-\gamma ), ~ \gamma \in B({\bf K},0,1) \} $ with $|\xi
_4|=p^{j+1}$ such that $F_j^1(\xi _1, \xi _2\gamma + \xi _3
(1-\gamma ))=((\xi _2+\xi _4)\gamma + (\xi _3+\xi _4)(1-\gamma ),0)$
for each $\gamma $, where $\xi _1,...,\xi _4\in {\bf K}$. Therefore,
one gets $\xi _2\gamma + \xi _3(1-\gamma )=0$ $\Leftrightarrow $ $(
\xi _2 - \xi _3) \gamma = - \xi _3$, $~ \gamma = \frac{\xi _3}{\xi
_3- \xi _2}$ hence $|\gamma |=1$ and $\gamma \in B({\bf K},0,1)$ and
there exists $z\in A_j$ with $|z_1|=|\xi _1|=p^{j-1}$ and $z_2=0$.
On the other hand, if $z\in B_j$, then $|z_1|=p^{j+1}$. Conditions
$(1,2)$ can be satisfied using a partition of $X$ into a disjoint
union of clopen balls.
\par We consider points $a\in {\bf K}_2\times ... \times {\bf K}_j$
and $b\in {\bf K}_1\times {\bf K}_{j+1}$ and $c=(c_1,c_2,...)\in X$
such that $c_1=0,$ ...,$c_j=0$ and put $\phi _j(a) = \xi \in \bf K$
when $j\ge 2$ such that $|\xi |=\max (0; q\in \Gamma _{\bf K}, ~
q\le 1- p^{j+1} \| a \| \} $, where $ \| a \| = \sup_{2\le l \le j}
|a_l|$, $a=(a_2,...,a_j)$, $a_l\in {\bf K}_l$ for each $l$. Put also
$\phi _1(a)=1$. Then a non-archimedean isotopy $H^j$ of $X$ onto $X$
exists so that $H^j_t(a,b,c) = (a, F^j_{t\phi _j(a)}(b),c).$ The
function $\phi _j(a)$ is continuous, since balls $B({\bf K}_2\times
... \times {\bf K}_j,x,r)$ are clopen in the normed space ${\bf
K}_2\times ... \times {\bf K}_j$ for each $x\in {\bf K}_2\times ...
\times {\bf K}_j$ and $r\in \Gamma _{\bf K}$.
\par The constructed sequence $H^j$ of non-archimedean isotopies
satisfies Conditions $1$ and $2$ of Lemma 19, hence \par $H^1_t = L
\prod_{j=0}^{\infty } H^{j+1} [1-p^j, 1-p^{j+1}]_t$ \\ is an
invertible non-archimedean isotopy pushing $x_0$ off $X$.
\par It remains to verify that $\{ H^j \} $ satisfies Conditions
19$(1,2)$. For this purpose take an arbitrary point $x\in X\setminus
\{ x_0 \} $ with $x=(0,...,0,a_j,a_{j+1},...)$, where $a_j$ is the
first non zero coordinate of $x$. For the composition $Q^j :=
H^j_1\circ ...H^2_1\circ H^1_1$ let $Q^j(x) =(b_1,...,b_j,
b_{j+1},a_{j+2},...)$ and
$Q^{j+1}(x)=(c_1,b_2,...,b_{j+1}c_{j+2},a_{j+3},...)$, where
$b_2=0$,...,$b_j=0$ for each $j\ge 2$. A neighborhood $U$ of $x$ and
an integer $k$ exist, when one of the coordinates
$b_2,...,b_{j+1},c_{j+2}$ is non zero, so that $\phi _k=0$ on $\pi
_{(2,...,k)}Q^k(U),$ where $\pi _{(l,...,k)}: X\to {\bf K}_l\times
... \times {\bf K}_k$ is a $\bf K$-linear projection with $l<k$.
Therefore, $H^v_t=id$ on $Q^k(U)$ for $v\ge k$. Particularly,
$Q^j(0,...,0,a_{j+2},...) = (\xi _4,0,...,0,a_{j+2},...)$, where
$|\xi _4| = p^{j+1}$. If $|b_1|\ne p^j$ and $c_{j+2}=0$, then
$|c_1|< p^{j+2}$ due to Condition $(2)$. In the case $|c_1|< p^{j+
2}$ a neighborhood $U$ of $x$ exists so that $|\pi _1\circ
Q^{j+1}(y)|<p^{j+2}$ for each $y\in U$. Thus $H^k_t=id$ on
$Q^{j+1}(U)$ for each $k\ge j+1$. That is, Condition 19$(1)$ is
fulfilled. Then Condition 19$(2)$ is satisfied, since
$H^{j+1}_t(y)=y$ when $|y_1|< p^j$, while $|\pi _1\circ Q^j(x_0)| =
p^{j+1}$.
\par {\bf 24. Lemma.} {\it Let $X=X_1\times X_2$, where either $X_1 ={\bf s}^{\alpha }$ and
$X_2={\bf s}^{\beta }$ or $X_1=c_0(\alpha ,{\bf K})$ and
$X_2=c_0(\beta ,{\bf K})$ with $card (\alpha )=card (\beta )=\aleph
_0$, $X_1$ and $X_2$ have origins $x_{0,1}$ and $x_{0,2}$
respectively. Suppose that $H$ is an invertible (non-archimedean)
isotopy pushing $x_{0,1}$ off $X_1$ and $\phi _r: X_2\to B({\bf
K},0,1)$ is a continuous one parameter family of maps with $r\in
B({\bf K},0,1)\setminus \{ 0 \} $ and $\phi _r^{-1}(1)=x_{0,2}$ and
$w: B({\bf K},0,1)\times B({\bf K},0,1)\to B({\bf K},0,1)$ is
continuous such that $w: B({\bf K},0,1)\setminus \{ 0 \} \times
B({\bf K},0,1)\setminus \{ 0 \} \to B({\bf K},0,1)\setminus \{ 0 \}$
is a mapping onto $B({\bf K},0,1)\setminus \{ 0 \}$ so that
$w^{-1}(1)=(1,1)$ and $w([\{ 0 \} \times B({\bf K},0,1)]\cup [B({\bf
K},0,1)\times \{ 0 \} ]) = \{ 0 \} $. Then $\mbox{}_rH_t(x,y) =
(H_{w(t,\phi _r(y))}(x),y)$ defines an invertible continuous one
parameter with $r\in B({\bf K},0,1)$ family of invertible
(non-archimedean) isotopies pushing the origin off $X$ for each
$r$.}
\par {\bf Proof.} In two considered cases $X$ is either $\bf s$ or
$c_0$. Since $\phi _r(y)$ is continuous on $(B({\bf K},0,1)\setminus
\{ 0 \} )\times X_2$, then the composite mapping $w(t,\phi _r(y)):
B({\bf K},0,1)\times [B({\bf K},0,1)\setminus \{ 0 \}]\times X_2 \to
B({\bf K},0,1)$ is continuous. But $H_b^{-1}$ is a continuous
(non-archimedean) isotopy. Therefore, $\mbox{}_rH_t^{-1}(x,y) =
(H^{-1}_{w(t,\phi _r(y))}(x),y) $ is the continuous non-arhimedean
isotopy. On the other hand, $\phi _r(x_{0,2})=1$ for each $r$,
consequently, $\mbox{}_rH_1(x_{0,1},x_{0,2}) =
(H_1(x_{0,1}),x_{0,2}) \ne x_0$ and $ \mbox{}_rH_1(X) = \{
H_{w(t,\phi _r(y))}(X_1),y): y \in X_2 \} = X$ for each $r$.
Moreover, $\mbox{}_rH_0(x,y) = (H_{w(0,\phi _r(y))}(x,y) =
(H_0(x),y) =(x,y) = id (x,y) ,$ since $w(0,b)=0$ for each $b\in
B({\bf K},0,1)$.
\par {\bf 25. Lemma.} {\it An invertible continuous non-archimedean
one parameter family of invertible isotopies $\mbox{}_r\tilde H$
exists with $r\in B({\bf K},0,1)\setminus \{ 0 \} $ each pushing the
origin $x_0$ off $\bf s$ so that $\mbox{}_r{\tilde H}_t$ is the
identity mapping outside a neighbourhood $U_r$ of $x_0$.}
\par {\bf Proof.} Let $\alpha $ and $\beta $ be two infinite subsets
in $\bf N$ such that $\beta = {\bf N}\setminus \alpha $. We define
the metric on $X^{\gamma } = {\bf s}^{\gamma }$ by the formula
$$d_{\gamma }(x,y) = \sum_{j\in \gamma } \min (p^{-j}, |x_j-y_j|)\in {\cal R},$$
where $x, y \in {\bf s}^{\gamma }$. Take in particular $\gamma
=\alpha $ or $\gamma = \beta $. Suppose wihout loss of generality
that $1, 2, 3\in \beta $, hence the diameter of ${\bf s}^{\alpha }$
is less than $p^{-3}$, $~diam ({\bf s}^{\alpha }) := \sup_{x,y\in
{\bf s}^{\alpha }} d_{\alpha }(x,y) <p^{-3}$. \par Next we consider
a non-archimedean isotopy $H$ pushing the origin off ${\bf
s}^{\alpha }$ (see Lemma 23). There exists a continuous map $\phi $
of ${\bf s}^{\beta }$ on $B({\bf K},0,1)$ so that $\phi
^{-1}(1)=x_{0,2}$ while $\phi (x)$ is zero for each $x\in {\bf
s}^{\beta }\setminus B({\bf s}^{\beta },x_{0,2},p^{-3})$, where
$B({\bf s}^{\beta },x_{0,2},q) = \{ y\in {\bf s}^{\beta }: ~
d_{\beta }(x,y) \le q \} $, $q\in \Gamma _{\bf K}$. This is
possible, since $B({\bf s}^{\beta },x_{0,2},q)$ is clopen in ${\bf
s}^{\beta }$ for $q>0$. Put $\mbox{}_1H_t(x,y) = (H_{w(t,\phi
(y))}(x),y)$ when $d((x,y),x_0)\ge p^{-2}$, since $1, 2, 3 \in \beta
$ and $d_{\beta }(y,x_{0,2})\ge p^{-3}$ and hence $\phi (y)=0$. We
denote by $l$ the least natural number in $\alpha $ and by $k$ a
natural number in $\beta $ greater than $l$. Let $\beta _1$ be the
subset of $\bf N$ formed from $\beta $ by the substitution $k\mapsto
l$, while $\alpha _1$ is made from $\alpha $ substituting $l$ with
$k$. \par There exists a family $F_{\lambda }$, with $\lambda \in
B({\bf K},0,1)$, of transfromations of $\bf s$ given by the formula:
$F_{\lambda }(x_1,x_2,...) = (y_1,y_2,...)$ with $y_j=x_j$ for each
$j\in {\bf N}\setminus \{ l, k \} $; whilst $y_l = (1-\lambda ) x_l$
when $|(1-\lambda ) x_l| > |\lambda x_k|$, $y_l = \lambda x_k$ when
$|\lambda x_k| \ge  |(1-\lambda ) x_l|$; $y_k=\lambda x_l$ when
$|\lambda x_l| \ge  |(1-\lambda ) x_k|$, $~ y_k = (1-\lambda ) x_k$
when $|(1-\lambda ) x_k| > |\lambda x_l|$.
\par Let $f(r)$ be a locally affine continuous mapping from $B({\bf
K},0,1)$ onto $B({\bf K},0,1)$ so that $f(B({\bf K},\xi _1,p^{-2}))=
B({\bf K},1,p^{-2})$ for $\xi _1\in \bf K$ with $|\xi _1|=p^{-1}$,
$f(B({\bf K},1,p^{-l-1})\setminus B({\bf K},1,p^{-l-2}))= B({\bf
K},0,p^{-l-1})\setminus B({\bf K},0,p^{-l-2})$ for each natural
number $l$, $ ~ f(1)=0$.
\par Define the mapping $\mbox{}_rH_t := F^{-1}_{f(r)}\circ \mbox{}_1H_t\circ F_{f(r)}$
for each $r\in W := [B({\bf K},\xi _1,p^{-2})\cup B({\bf
K},1,p^{-2})]$. If $d(x,x_0)\ge p^{-2}$ and $r\in W$, then
$d(F_{f(r)}(x),x_0)\ge d(x,x_0)$ and hence $\mbox{}_rH_t(x)=x$,
since $|1-\lambda |=1$ for any $\lambda \in B({\bf K},0,p^{-1})$ and
$|\lambda |=1$ for each $\lambda \in B({\bf K},1,p^{-1})$. \par Now
we define $\mbox{}_rH$ for $r\in B({\bf K},1,p^{-1})\setminus W =:
C_1$ so that $\mbox{}_rH_t(x,y) = (S_{w(t,\phi _r(y))}(x),y)$ for
each $x\in {\bf s}^{\alpha _1}$ and $y\in  {\bf s}^{\beta _1}$,
where $T$ is a linear over $\bf K$ operator on $\bf s$ interchanging
a finite number of coordinates, $S_t := T^{-1}\circ H_t\circ T$.
\par Let $\xi _1$ be a marked point as above
and define $\phi _{\xi _1}$ as a map of ${\bf s}^{\beta _1}$ onto
$B({\bf K},0,1)$ such that $\phi _{\xi _1}^{-1}(1)=x_{0,2}$ and
$\phi _{\xi _1}(y)=0$ for each $y\in {\bf s}^{\beta _1}\setminus
U_1$, where $U_1$ is a small neighborhood of $x_{0,2}$, $~U_1= \{
y\in {\bf s}^{\beta _1}: ~ d_{\beta _1} (x_{0,2},y) \le p^{-3} \} $.
This mapping $\phi _{\xi _1}$ can be presented as the composition
$\phi _{\xi _1} = \phi \circ T_1$, where a mapping $T_1 : {\bf
s}^{\beta _1}\to {\bf s}^{\beta _1}$ is given by the formula
$T_1(x_1,x_2,..) = (y_1,y_2,...)$ with $y_j=x_j$ for each $j\notin
\{ l, k \} $, also $y_l=x_k$ and $y_k=-x_l$.
\par If $d(z,x_0)\ge p^{-2}$, $~z=(x,y)$, $~ x\in {\bf s}^{\alpha
_1}$, $~y\in {\bf s}^{\beta _1}$, then $\phi _{\xi _1} (y)=0$, since
$d(y,x_{0,2})\ge p^{-2}$ and $d(T_1(y),x_{0,2})\ge p^{-2}$. For $\xi
_v\in B({\bf K},0,p^{-v}) \setminus B({\bf K},0,p^{-v-1}) =: C_v$
with a natural number $2\le v\in {\bf N}$, let $\phi _v $ be a
continuous mapping of ${\bf s}^{\beta _1}$ onto $B({\bf K},0,1)$ so
that $|\phi _{\xi _{v+1}}(y)| \le |\phi _{\xi _{v}}(y)|$ for each
$y\in {\bf s}^{\beta _1}$, also $\phi _{\xi _v}^{-1}(1) = x_{0,2}$,
$ ~ \phi _{\xi _v}=0$ for each $y$ with $d(y,x_{0,v})>p^{-v-2}$,
where $x_{0,v}$ is the origin in ${\bf s}^{\beta _v}$, $~T_v$ and
$\alpha _v$ and $\beta _v$ are defined by induction. For each $y\in
{\bf s}^{\beta _v}$ let $$\frac{\phi _{\xi _v}(y)-\phi _r(y)}{\phi
_{\xi _{v-1}}(y) - \phi _r(y)} = \frac{\xi _v - r}{\xi _{v-1}-r } $$
when $\phi _{\xi _v}(y) \ne  \phi _{\xi _{v-1}}(y)$ for each $r\in
C_v\setminus \{ \xi _v \} $, while $\phi _r(y) = \phi _{\xi _v}(y)$
when $\phi _{\xi _v} (y) =\phi _{\xi _{v-1}}(y)$. Put
$\mbox{}_rH_t(x,y) = (S_{w(t,\phi _r(y))}(x),y)$ for each $r\in C_v$
by induction on $2\le v \in {\bf N}$.
\par On the other hand, due to Lemmas 6 and 15, Theorem 16 there
exists a homeomorphism $\theta $ of topological spaces from $B({\bf
K},0,1)\setminus \{ 0 \} $ onto $[B({\bf K},\xi _1,p^{-2})\cup
B({\bf K},0,p^{-2})\cup B({\bf K},1,p^{-2})]\setminus \{ 0 \} $ such
that $\theta (\xi _1) = \xi _1$, $~ \theta (1)=1$ and $\lim_{t\to 0}
\theta (t)=0$, hence $\mbox{}_{\theta (r)}H_t$ is the claimed
isotopy $\mbox{}_r{\tilde H}_t$.
\par {\bf 26. Lemma.} {\it There exists an invertible non-archimedean isotopy
$F$ pushing a point off $A_0$.}
\par {\bf Proof.} We consider the set $c_0(1) := \{ x\in c_0: ~
x_1=1, ~ x=(x_1,x_2,...), ~\forall j \in {\bf N} ~ x_j \in {\bf K}
\} $ and points $q_j = (q_{j,1},...,q_{j,j},0,0,..)\in c_0$,
$q_{j,1}=1$,...,$q_{j,j}=1$ for each $j\in {\bf N}$. Take the
neighborhood $U_j = \{ (1,x_2,...)\in c_0(1): ~ \max_{i=2}^j |1-x_i|
< p^{-j} \} $ of the point $q_j$. Define $H^1_t(x)=x+te_2$,
$H^i_t(x) = x + p^{1-i} t \eta _i(x) e_{i+1}$, where $\eta _i(x)\in
\bf K$, $~|\eta _i(x)|=d(x,c_0(1)\setminus U_i)$ for each $i\ge 2$,
and \par
$H_t= L \prod_ {i=0}^{\infty } H^{i+1} [1-p^i,1-p^{i+1}]_t$, \\
where $e_i=(0,...,0,1,0,...)$ is the vector with unit $i$-th
coordinate and zero others. In view of Theorem 4.1 \cite{anderbing}
and Lemma 19 $H$ is an invertible non-archimedean isotopy pushing
$q_1$ out of $c_0(1)$. \par The topological space $c_0(1)$ is the
closed subset in $c_0$ and is a union of $\omega ({\bf K})|\Gamma
_{\bf K}| \aleph _0$ disjoint balls $B(c_0(1),x,r)$, where $x\in
c_0(1)$, $~r\in \Gamma _{\bf K}$, $~\omega (P)$ denotes the
topological weight of a topological space $P$, whilst $card (S) =
|S|$ denotes the cardinality of a set $S$. Then the topological
space $A_0$ is also the closed subset in $c_0$. The topological
spaces $c_0$ and $c_0\oplus {\bf K}$ are linearly homeomorphic, but
$c_0\oplus {\bf K}$ is also linearly homeomorphic with $c$,
consequently, $c_0$ and $c$ are linearly homeomorphic. Consider the
homemorphism $\nu : c_0\to c$ of $c_0$ onto $c$.
\par Particularly, consider linear topological isomorphism
$\nu : c_0\to c$ such that $\nu (x)=y$ with $y_k = \sum_{j=1}^k x_j$
for each $k\in \bf N$, since the series $\sum_{j=1}^k x_j$ converges
if and only if $\lim_j x_j=0$ due to the ultrametric inequality.
Moreover, the topological space $\nu (A_0)=: W_0$ is the disjoint
union of $\omega ({\bf K}) |\Gamma _{\bf K}| \aleph _0$ balls
$B(c,y,r)$ with $y\in W_0$ and $\Gamma _{\bf K}\ni r\le p^{-1}$,
since from $y\in W_0$ and $z\in B(c,0,r)$ it follows $|y_k+z_k|\le
\max (|y_k|, |z_k|)\le 1$ for each natural number $k$ and from
$|y_l|=1$ it follows $|y_l +z_l|=1$, while $\nu (A_0) = \{ y\in c: ~
\sup_k |y_k|=1 \} $. \par Therefore, each two balls $B(c_0(1),x,r)$
and $B(c,y,q)$ with $r, q\in \Gamma _{\bf K}$ are homeomorphic. Thus
$A_0$ and $c_0(1)$ are homeomorphic. Therefore, there exists an
invertible non-archimedean isotopy $F_t(y) = \psi ^{-1}\circ H_t
\circ \psi (y)$, where $\psi $ is a homeomorphism of $A_0$ onto
$c_0(1)$ described above.
\par {\bf 27. Corollary.} {\it There exists an invertible
non-archimedean isotopy $G$ pushing a point off $c_0$.}
\par {\bf Proof.} The topological spaces $c_0$ and $c_0(1)$ are
homeomorphic with a homeomorphism $\eta : c_0\to c_0(1)$. Indeed,
the topological space $c_0$ can be presented as the disjoint union
of $\omega ({\bf K}) |\Gamma _{\bf K}|\aleph _0$ balls $B(c_0,x,r)$
with $x_0\in c_0$ and $r\in \Gamma _{\bf K}$. From \S 26 it follows
that $G_t(y) = \eta ^{-1}\circ H_t\circ \eta (y)$ is an invertible
non-archimedean isotopy pushing a point off $c_0$.
\par {\bf 28. Lemma.} {\it There exists an invertible non-archimedean
isotopy $H$ pushing the origin $x_0=0$ off $c_0$ so that $H_t$ is
the identity outside the unit ball in $c_0$ containing zero for each
$t\in B({\bf K},0,1)$.}
\par {\bf Proof.} The unit ball $B(c_0,0,1)$ is clopen in $c_0$. On
the other hand, $B(c_0,0,1)$ is homeomorphic with $A_0$. An
invertible non-archimedean isotopy evidently has an extension from
$B(c_0,0,1)$ onto $c_0$. From Lemma 26 and the equality
$B(c_0,0,1)=B(c_0,x,1)$ for each $x\in B(c_0,0,1)$, particularly,
for $\| x \| =1$, the statement of this lemma follows.
\par {\bf 29. Lemma.} {\it There exists an invertible continuous
one parameter $r\in B({\bf K},0,1)$ family of invertible
non-archimedean isotopies $\mbox{}_rH$ each pushing the origin $x_0$
off $c_0$ such that $\mbox{}_rH_t$ is the identity outside
$B(c_0,0,|r|)$ for each $t\in B({\bf K},0,1)$.}
\par {\bf Proof.} If $r\in {\bf K}\setminus \{ 0 \} $, then $m_r:
c_0\to c_0$ is the linear homeomorphism of $c_0$ onto $c_0$, where
$m_r(x)=rx$ for any $x\in c_0$. Therefore, the desired isotopy is
given by the formula: $\mbox{}_rH_t=m_r\circ H_t\circ m_{1/r}$,
where $H_t$ is the isotopy provided by Lemma 28.
\par {\bf 30. Lemma.} {\it Let $X$ be one of the topological spaces $c_0$
or $\bf s$, let also $\Omega $ be an open covering of $X$. Suppose
that $\alpha $ is an infinite proper subset of $\bf N$ and $Q$ is a
closed subset in $X$ deficient with respect to $\alpha $. Then for
any open subset $U$ containing $Q$ a homeomorphism $g$ of
$X\setminus Q$ onto $X$ exists such that $g$ is limited by $\Omega
$, $~ h|_{X\setminus U} = id$, $~\pi _j(x) = \pi _j(g(x))$ for each
$j\in \beta $, where $\beta = {\bf N}\setminus \alpha $.}
\par {\bf Proof.} In view of Lemma 4 there exists the decomposition
of $X$ into the product of two topological spaces $X=X^{\alpha
}\times X^{\beta }$, where either $X^{\alpha } = c_0(\alpha ,{\bf
K})$ or $X^{\alpha } = {\bf s}^{\alpha }$ for either $X=c_0$ or
$X={\bf s}$ respectively. \par Without loss of generality we
consider the case $\pi _{\alpha }(0)=x_{0,\alpha }$, where
$x_{0,\alpha } =0 \in X^{\alpha }$. Put $\pi _{\beta }(Q)=Q'$, hence
$Q'\subset X^{\beta }$.
\par Take an open covering $\cal W$ of $X$ of mesh less than one
relative to the norm on $c_0$ or the metric $d$ on $\bf s$ so that
\par $(1)$ $\cal W$ is a refinement of $\Omega $,
\par $(2)$ if $S\in \cal W$ and $S\cap Q\ne \emptyset $, then
$S\subset U$. \par Then let $v(y)$ be a function defined by the
formula:
\par $(3)$ $v(y)= \sup \{ \epsilon: ~ \epsilon \in \Gamma _{\bf K},
~ \exists ~ S\in {\cal W} ~ B(X,(x_{0,\alpha },y),p \epsilon
)\subset S \} $, where $y\in X^{\beta }$. This function satisfies
the inequality $v(y)>0$ for each $y\in X^{\beta }$.
\par For a set $A= \{ y \in X^{\beta }: ~ B(X,(x_{0,\alpha
},y),v(y))$ $\mbox{ is not contained in }$ $U \} $ we define the
function
\par $(4)$ $\tau (y) := d(y,A)/(d(y,A) + d(y,Q'))$, when
$A$ is non void, $A\ne \emptyset $, while $\tau (y)=1$ for
$A=\emptyset $. This definition implies that $(cl ~ A)\cap
Q'=\emptyset $, $\tau (y)=1$ for each $y\in Q'$ and $\tau (y)=0$ for
any $y\in A$.
\par In accordance with Lemmas 25 and 29 an invertibly continuous
family $\mbox{}_rH$ with $r\in B({\bf K},0,1)\setminus \{ 0 \} $ of
invertible non-archimedean isotopies exist each pushing $x_{0,\alpha
}$ off $X^{\alpha }$ and so that $\mbox{}_rH_t|_{X^{\alpha
}\setminus B(X^{\alpha },x_{0,\alpha },|r|)}=id$.
\par We then define the mapping
\par $(5)$ $g(x,y) = (\mbox{}_{w(y)}H_{\mu (y)}(x),y)$ and verify below
that it satisfies the desired properties, where $w(y)$ and $\mu (y)$
are continuous mappings from $X^{\beta }$ into $\bf K$ such that
$\frac{v(y)}{p} \le |w(y)| \le v(y)$ and $\frac{\tau (y)}{p} \le
|\mu (y)| \le \tau (y)$ for each $y\in X^{\beta }$, where $\mu
(y)=1$ for each $y\in Q'$. This implies that $\mbox{}_{w(y)}H_{\mu
(y)}$ is a homeomorphism of $(X^{\alpha }\setminus \{ x_{0,\alpha }
\} )\times \{ y \} $ onto $X^{\alpha }\times \{ y \} $, when $y\in
Q'$, since $\mu (y)=1$ for each $y\in Q'$. Then
$\mbox{}_{w(y)}H_{\mu (y)}$ is a homeomorphism of $X^{\alpha }\times
\{ y \} $ onto $X^{\alpha }\times \{ y \} $, since $\tau (y)<1$ for
any $y\in X^{\beta }\setminus Q'$.
\par As the composition of continuous mappings, the mapping $g(x,y)$
given by Formula $(5)$ also is continuous, since $\mbox{}_rH$ is a
continuous family of isotopies. The inverse mapping is given by the
formula $g^{-1}(x,y) = (\mbox{}_{w(y)}H^{-1}_{\mu (y)} (x),y).$
Since $\mbox{}_rH$ is an invertibly continuous family of invertible
non-archimedean isotopies, then $g^{-1}$ is continuous. Thus $g$ is
a homeomorphism of $X\setminus Q$ onto $X$.
\par Formula $(5)$ implies that $\pi _j(z)=\pi _j(g(z))$ for every
$z\in X$ and $j\in \beta $, since $z=(x,y)$ with $x\in X^{\alpha }$
and $y\in X^{\beta }$.
\par On the other hand, the mapping $g$ is the identity on
$X\setminus U$, since
\par $(6)$ $\mbox{}_rH_t|_{X^{\alpha }\setminus B(X^{\alpha
},x_{0,\alpha },|r|)} = id $ and $B(X,(x_{0,\alpha },y),|w(y)|)
\subset U$ when $\tau (y)\ne 0$. Moreover, Condition $(6)$ implies
that either $g(x,y) = (x,y)$ or $x$ and $\mbox{}_{w(y)}H^{-1}_{\mu
(y)}(x)\in B(X^{\alpha },x_{0,\alpha },|w(y)|)$. But the definition
of $w(y)$ means that $(x,y)$ and $g(x,y)$ belong to an element $S$
of a covering $\nu $ containing $(x_0,y)$, i.e. $(x_0,y)\in S$.
\par {\bf 31. Theorem.} {\it Let $X$ be either $c_0$ or $\bf s$, let also
$U$ be an open subset in $X$. Suppose that $ \{ K_j: ~ j\in {\bf N}
\} $ is a sequence of closed subsets of $X$ so that $K_j\subset U$
and $K_j$ has an infinite deficiency for each $j$. Then a
homeomorphism $g$ of $X\setminus \bigcup_{j=1}^{\infty }K_j$ onto
$X$ exists so that $g|_{X\setminus U} = id$.}
\par {\bf Proof.} The definition of the infinite deficiency implies
that for each $K_j$ an infinite subset $\beta _j\subset \bf N$
exists so that $K_j$ is deficient with respect to $\beta _j$. The
family $ \{ K_j: j \in {\bf N} \} $ is countable and $\pi _l(K_j)$
consists of a single element for each $l\in \beta _j$, consequently,
there exists a disjoint family of infinite subsets $\alpha _j\subset
\beta _j$, $~\alpha _i\cap \alpha _j= \emptyset $ for each $i\ne j$,
such that $K_j$ is infinite deficient with respect to $\alpha _j$.
\par In accordance with Lemma 30 a sequence of homeomorphisms
$g_j$ and a sequence of coverings $G_j$ with $j\in \bf N$ satisfying
conditions of Theorem 4.3 \cite{anderbing} exist. Mention that
results of \S 4 \cite{anderbing} are also valid for a metric space
with a metric $d$ satisfying the non-archimedean inequality and with
values in ${\cal R}$ substituting in proofs $1/2$ on $1/p<1\in {\cal
R}$, since a ring $\cal R$ is complete as a uniform space. That is,
for each natural number $j\ge 1$ there exists a homeomorphism $g_j$
of $X\setminus g_{j-1}\circ ... \circ g_1\circ id (K_j\setminus
\bigcup_{l=1}^{j-1}K_l)$ onto $X$ with $g_0=id$ such that $g_j(x)=y$
for each $x\in X\setminus g_{j-1}\circ ... \circ g_1\circ id
(K_j\setminus \bigcup_{l=1}^{j-1}K_l)$ with $x_k=y_k$ for any $k\in
{\bf N}\setminus \alpha _j$. Therefore, the mapping
$L\prod_{j=1}^{\infty } g_j$ is a homeomorphism of $X\setminus
\bigcup_{j=1}^{\infty }K_j$ onto $X$.
\par {\bf 32. Corollary.} {\it The topological space $c_0\setminus
\bigcup_{j=1}^{\infty }E^j$ is homeomorphic with $c_0$.}
\par {\bf 33. Corollary.} {\it If $ \{ C_j: ~ j\in {\bf N} \} $
is a countable family of compact subsets of $\bf s$, then ${\bf
s}\setminus \bigcup_{j=1}^{\infty } C_j$ is homeomorphic to $\bf
s$.}
\par {\bf Proof.} A homeomorphism $g$ of $\bf s$ onto itself so that
$g(C_j)$ is infinitely deficient for each $j\in \bf N$ exists due to
Theorem 21. Then from Theorem 31 the assertion of this corollary
follows.
\par {\bf 34. Corollary.} {\it Let $U$ be an open subset of $c_0$,
let also $K_j$ be a compact subset in $U$ for each natural number
$j$. Then a homeomorphism $g$ of $c_0\setminus \bigcup_{j=1}^{\infty
} K_j$ onto $c_0$ exists so that $g|_{c_0\setminus U} = id$.}
\par {\bf Proof.} In view of Lemma 20 and Theorem 21 a homeomorphism
$\eta $ of $c_0$ onto $c_0$ exists such that $\eta (K_j)$ is
infinitely deficient. Then by Theorem 31 there exists a
homeomorphism $\xi $ of $c_0\setminus \bigcup_{j=1}^{\infty } g
(K_j)$ onto $c_0$ so that $\xi |_{c_0\setminus \eta (U)} = id$,
consequently, $g= \eta ^{-1}\circ \xi \circ \eta $.
\par {\bf 35. Lemma.} {\it There exists a homeomorphism
of $c_0\setminus \bigcup_{n=1}^{\infty } E^n$ with $A_2^*$.}
\par {\bf Proof.} The topological space
$A_2 := \{ x\in c_0: ~ \sup_{1\le j\le k\in {\bf N}}
|\sum_{j=1}^kx_j| =1, ~ \sum_{j=1}^kx_j\ne 1 ~ \forall k\in {\bf N},
~ \sum_{j=1}^{\infty }x_j=1 \} $ is homeomorphic with $A_3 := \{
y\in c: ~ \sup_{k\in {\bf N}} |y_k| =1, y_k\ne 0 \forall k \} $.
Indeed, the mapping $h_1(x)=y$ such that $y_k= (\sum_{j=1}^kx_j)$ is
the required homeomorphism, since $x_1=1-\sum_{j=2}^{\infty }x_j$
and $|x_1|\le |1|$, $\lim_k y_k = (\sum_{j=1}^{\infty } x_j) =1$.
Therefore, the space $A_2^*$ is homeomorphic with $A_3^* =
A_3\setminus \bigcup_n E^n$, since $x_n=0$ is equivalent to
$y_n=y_{n-1}$ and $E^n=E^1\oplus E^{n-1}$, where $E^n$ is
homeomorphic with ${\bf K}^n$.
\par The topological spaces $c$ and $A_3$ are homeomorphic as
disjoint unions of $|\Gamma _{\bf K}| \aleph _0 \omega ({\bf K})$
balls $B(c,x,q)$ with $q\le 1/p$, $~q\in \Gamma _{\bf K}$, where
$x\in c$ or $x\in A_3$ respectively. On the other hand, $c_0$ and
$c$ are homeomorphic (see \S 26), whilst the topological spaces
$c_0$ and $c_0\setminus \bigcup_{j=1}^{\infty }E^j$ are homeomorphic
by Corollary 32. Thus $c_0\setminus \bigcup_{j=1}^{\infty }E^j$ is
homeomorphic with $A_3^*$ due to Theorem 31 and hence with $A_2^*$.
\par {\bf 36. Lemma.} {\it The topological spaces $s_*$ and
$s$ are homeomorphic.}
\par {\bf Proof.} Consider the topological space $s$ realized as in \S 7.
The mapping $h(x)=y$ with $y_k=x_k-1$ for each $k\in \bf N$ is the
isomorphism of $s$ with $$s_1 := \prod_{j=1}^{\infty } (B({\bf
K},0,1)_j\setminus \{ 0 \} ).$$ Put $C_1^j := \{ y\in P: ~ y_k=-1
\mbox{ or } y_k=0 ~ \forall k>j \} $, where $$P=\prod_{j=1}^{\infty
} B({\bf K},0,1)_j.$$ Since a field $\bf K$ is locally compact and
the unit ball $B({\bf K},0,1)$ of it is compact, then the
topological space $P$ is compact and each subset $C_1^j$ in it is
compact. In view of Theorem 4.3 \cite{anderbing}, Remark 7 and Lemma
25 and Corollary 33 the topological spaces $s_1$ and $s_{1*} := s_1
\setminus \bigcup_{j=1}^{\infty } C_1^j$ are homeomorphic.
Therefore, the topological spaces $s_*$ and $s$ are homeomorphic.
\par {\bf 37. Remark.} Thus Theorem 2 follows from the sequence of
homeomorphisms $c_0 \approx c_0\setminus \bigcup_jE^j \approx A_2^*
\approx A_1^*\approx s_* \approx s$ demonstrated above in Corollary
32, Lemmas 35, 12, 10 and 36.
\par Theorem 2 has a generalization over a field $\bf F$ of
separable type over $\bf K$, i.e. when $\bf F$ has a multiplicative
norm extending that of $\bf K$ such that $\Gamma _{\bf K}\subset
\Gamma _{\bf F}$ and $\bf F$ has an equivalent norm $|*|_{1,\bf F}$
not necessarily multiplicative so that $|a|_{\bf F}/p \le |a|_{1,\bf
F} \le p |a|_{\bf F}$ and $|a|_{1,\bf F} \in \Gamma _{\bf K}\bigcup
\{ 0 \} $ for each $a\in \bf F$ and $({\bf F}, |*|_{1,\bf F})$ is
isomorphic as the Banach space over $\bf K$ with $c_0(\omega _0,{\bf
K})$. Indeed, $c_0(\omega _0,c_0(\omega _0,{\bf K}))$ is isomorphic
with $c_0(\omega _0,{\bf F})$ and $({\bf K}^{\omega _0})^{\omega
_0}$ is isomorphic with $(c_0(\omega _0,{\bf K}))^{\omega _0}$ and
hence with ${\bf F}^{\omega _0}$, consequently, the topological
spaces $c_0(\omega _0,{\bf F})$ and ${\bf F}^{\omega _0}$ are
homeomorphic.

\par {\bf 38. Corollary.} {\it If the topological $\bf K$ linear space
$c_0(\alpha )$ with $card (\alpha ) >\aleph _0$ is supplied with the
projective limit topology $\tau _{pr}$ induced by ${\bf K}$ linear
projection operators $\pi ^{\alpha }_{\beta }: c_0(\alpha )\to
c_0(\beta )$ associated with the standard base $\{ e_j: ~ j\in
\alpha \} $ for each $\beta \subset \alpha $ with $card (\beta ) =
\aleph _0$, then $(c_0(\alpha ), \tau _{pr})$ is topologically
homeomorphic with $({\bf K}^{\alpha }, \tau _{ty})$.}
\par {\bf Proof.} Mention the fact that the topological space
$({\bf K}^{\alpha },\tau _{ty})$ can be presented as the projective
limit of topological spaces $({\bf K}^{\beta }, \tau _{ty})$ with
$\beta \subset \alpha $, $card (\beta )=\aleph _0$, when $card
(\alpha )>\aleph _0$ (see also \cite{eng}). Therefore, applying
Theorem 2 we get this corollary.
\par {\bf 39. Conclusion.} The results of this paper can be used for subsequent studies of
non-archimedean functional analysis, measure theory, topological
algebra and geometry, topological vector spaces and manifolds on
them.

\end{document}